\documentclass[12pt]{amsart}
\usepackage[utf8]{inputenc}
\usepackage{graphicx}
\usepackage{amsfonts}
\usepackage[portuguese,english]{babel}
\usepackage{amscd,color}
\usepackage{amsmath}
\usepackage{amssymb}
\usepackage{enumerate}

\usepackage[
    colorlinks=true,       
    linkcolor=blue,        
    citecolor=blue,        
    urlcolor=blue,         
    bookmarksnumbered=true,
    bookmarksopen=true,
    bookmarksopenlevel=3,
    pdfstartview=FitH,
    pdfmenubar=true,
    pdftoolbar=true,
    bookmarks=true,
    filecolor=magenta,
    plainpages=false,
    pdfpagelabels,
    breaklinks
]{hyperref}

\usepackage{float}
\usepackage{fancyhdr}
\usepackage[
  lmargin=3.0cm,
  rmargin=3.0cm,
  tmargin=3.0cm,
  bmargin=3.0cm
]{geometry}

\usepackage{amsthm}
\usepackage{xcolor}

\makeatletter
\@namedef{subjclassname@2010}{%
	\textup{2010} Mathematics Subject Classification}
\makeatother

\theoremstyle{plain}

\newtheorem{corollary}{\bf Corollary}

\newtheorem{lemma}{\bf Lemma}

\newtheorem{proposition}{\bf Proposition}
\newtheorem{remark}{Remark}

\newtheorem{theorem}{\bf Theorem}
\numberwithin{equation}{section}

\newcommand{\blue}[1]{{\color{blue}#1}}

\newtheorem{maintheorem}{\bf Theorem}

\newtheorem{maincorollary}{\bf Corollary}

\title[Integral formulas and Serrin-type problems for the drift Laplacian]{Integral formulas and Serrin-type problems for the drift Laplacian}

\author[M. Andrade]{Maria Andrade}
\address[M. Andrade]{Universidade Federal de Sergipe - UFS, Departamento de Matemática, 49100-000, S\~ao Cristov\~ao-SE, Brazil.
}
\email{\href{mailto: maria@mat.ufs.br}{ maria@mat.ufs.br}}

\author[A. Freitas]{Allan Freitas}
\address[A. Freitas]{Federal University of Paraíba, Department of Mathematics, 58051-900, João Pessoa-PB, Brazil}

\address{Current address: University of São Paulo, Institute of Mathematics, Statistics and Computer Science,  05508-090, São Paulo-SP, Brazil}

\email{\href{mailto: allan@mat.ufpb.br}{allan@mat.ufpb.br}; \href{mailto: allan.freitas@academico.ufpb.br}{allan.freitas@academico.ufpb.br}; \href{mailto: allan.freitas@ime.usp.br}{allan.freitas@ime.usp.br}}

\author[J. S. Sindeaux]{Joyce S. Sindeaux}
\address[J. S. Sindeaux]{Regional University of Cariri, 63105-000, Crato – CE, Brazil}
\email{\href{mailto: }{joyce.sindeaux@urca.br}}


\begin{document}
\begin{abstract}
We study boundary value problems for the drift Laplacian on weighted manifolds. For Euclidean domains with a positive homogeneous weight, we establish weighted integral formulas for the torsion problem and use them to prove weighted versions of Serrin’s symmetry theorem and Alexandrov’s soap bubble theorem under curvature assumptions of integral type. On general weighted manifolds, we obtain estimates for the weighted torsion function and show that, unlike the unweighted case, they are strict for smooth positive weights. Finally, we prove a general integral identity for solutions of the weighted Poisson equation and apply it to a reverse Serrin problem.
\end{abstract}
\maketitle


\section{Introduction}

Let $\Omega\subset\mathbb{R}^n$ be a bounded domain with smooth boundary and let $u$ be its
\emph{torsion function}, that is, the solution of
\begin{equation}\label{intro:torsion}
    \Delta u=-1\quad\text{in }\Omega,\qquad u=0\quad\text{on }\partial\Omega .
\end{equation}
The quantity $T(\Omega)=\int_\Omega u\,dx$ is called the \emph{torsional rigidity} of
$\Omega$. Both have a clear physical meaning: in elasticity, $T(\Omega)$ measures the
resistance to twisting of a beam with cross section $\Omega$, and in fluid mechanics $u$ is
the velocity profile of a viscous incompressible fluid flowing along a straight pipe; see
\cite{PS,soko,sperb}.

In a celebrated paper, Serrin~\cite{Serrin1} proved that if the solution of
\eqref{intro:torsion} also satisfies the overdetermined condition
\begin{equation}\label{intro:neumann}
    u_\nu=-c\quad\text{on }\partial\Omega,
\end{equation}
where $\nu$ is the outward unit normal and $c$ is a constant, then $\Omega$ is a ball and $u$
is radially symmetric. Serrin's proof is based on
the method of moving planes, introduced by Alexandrov~\cite{Alexandrov2} to prove his Soap
Bubble Theorem: \emph{spheres are the only closed embedded hypersurfaces of constant mean
curvature in Euclidean space}. Weinberger~\cite{Weinberger}
gave a short proof of Serrin's theorem, which combines the maximum principle for the
$P$-function $P=|\nabla u|^2+\frac{2}{n}u$ with the Pohozaev identity~\cite{pohozaev}.
In the same spirit, Reilly~\cite{reilly} and Ros~\cite{ros} obtained proofs of Alexandrov's
theorem based on integral formulas. We refer to \cite{Magnanini2017,NT} for surveys on the
relationship between these two symmetry results.

The integral approach has proved to be remarkably flexible. On the one hand, it does not rely
on the invariance of the operator under reflections, and therefore it has been the main tool
for extending Serrin's theorem to space forms, warped products and manifolds carrying a
conformal vector field \cite{CV,FR_2022,FAM_2023,Roncoroni2019}, as well as to convex cones
\cite{PacTra,CirRon}. On the other hand, it produces \emph{quantitative} information.
Magnanini and Poggesi \cite{Magnanini2019,Magnanini2020} showed that Weinberger's and Reilly's
arguments can be condensed into two integral identities. With the normalization
\eqref{intro:torsion}, the first one reads
\begin{equation}\label{intro:MP}
    -2\int_\Omega u\,|\mathring{\nabla}^2u|^2\,dx
    =\int_{\partial\Omega}\big(u_\nu^2-c^2\big)\big(u_\nu-q_\nu\big)\,d\sigma,
\end{equation}
where $\mathring{\nabla}^2u=\nabla^2u-\frac{\Delta u}{n}\,g$ is the traceless Hessian of $u$,
$c$ is any constant and $q$ is a quadratic polynomial with $\Delta q=-1$. The second one relates
the same volume term to the deviation of the mean curvature of $\partial\Omega$ from a
reference constant. The quantity $|\mathring{\nabla}^2u|^2$, the \emph{spherical} (or
Cauchy--Schwarz) \emph{deficit}, vanishes identically if and only if $u$ is a quadratic
polynomial with spherical level sets. Hence \eqref{intro:MP} gives Serrin's theorem at once
when \eqref{intro:neumann} holds. More importantly, the right-hand side of \eqref{intro:MP}
measures how far $u_\nu$ is from being constant, and this is the starting point of the sharp
stability estimates for Serrin's problem and for the Soap Bubble Theorem obtained in
\cite{Magnanini2019,Magnanini2020,MPnearly}; see also \cite{ABR,BNST,CM} for earlier approaches
to stability, and \cite{poggesi} for convex cones. A Riemannian counterpart of
\eqref{intro:MP}, on manifolds endowed with a conformal vector field, was recently obtained
in \cite{AF}.

In this paper we investigate these questions for the \emph{drift Laplacian}
\begin{equation*}
    \Delta_wu=\Delta u+\langle\nabla\log w,\nabla u\rangle=\frac1w\operatorname{div}(w\nabla u),
\end{equation*}
where $w$ is a positive function. This is the natural Laplacian of the weighted manifold
$(M^n,g,d\mu_w=w\,d\mathrm{vol})$: it is symmetric with respect to $d\mu_w$, and its geometry
is governed by the $\alpha$-Bakry--\'Emery--Ricci tensor
\begin{equation*}
    \operatorname{Ric}^\alpha_w=\operatorname{Ric}-\nabla^2\log w
    -\frac1\alpha\,d\log w\otimes d\log w,\qquad\alpha>0,
\end{equation*}
through a Bochner formula in which $n+\alpha$ plays the role of the dimension
\cite{bakry,Lott_2003,wey}. Weighted manifolds (or manifolds with density) provide a unified
framework for isoperimetric and concentration inequalities, diffusion processes and optimal
transport; see \cite{morgan,gromov,villani}. There is also a direct physical motivation. Since
the weighted torsion problem
\begin{equation}\label{intro:problem}
    \Delta_wu=-1\quad\text{in }\Omega,\qquad u=0\quad\text{on }\partial\Omega
\end{equation}
is equivalent to $\operatorname{div}(w\nabla u)=-w$, it describes the analogue of
\eqref{intro:torsion} in an inhomogeneous medium whose conductivity (or permittivity,
in the electrostatic interpretation; see \cite{jackson}) is the function $w$.
The associated \emph{weighted torsional rigidity} is
\begin{equation*}
    T_w(\Omega)=\int_\Omega u\,d\mu_w=\int_\Omega|\nabla u|^2\,d\mu_w .
\end{equation*}

We mainly consider the Euclidean space endowed with a weight that is positive and \emph{homogeneous of degree
$\alpha>0$}, that is, $\langle\nabla w,x\rangle=\alpha w$, where $x$ is the position vector;
the model example is $w=|x|^\alpha$. Homogeneous weights appear naturally in the sharp
weighted isoperimetric inequalities in convex cones of Cabr\'e, Ros-Oton and
Serra~\cite{Cabre} and in the classification of stable hypersurfaces
of Ca\~nete and Rosales~\cite{CC_2014}. For such weights 
\begin{equation*}
    q(x)=-\frac{|x|^2}{2(n+\alpha)}
\end{equation*}
solves $\Delta_wq=-1$, and balls centered at the origin support solutions of
\eqref{intro:problem} with constant normal derivative. The converse is the content of the
theorem of Ruan, Wang and Huang~\cite{RWH_2018}: if \eqref{intro:problem} admits a solution
with $u_\nu$ constant on $\partial\Omega$ and such that
$\operatorname{Ric}^\alpha_w(\nabla u,\nabla u)\geq0$ \emph{pointwise} in $\Omega$ and
$\nabla u(0)=0$, then $\Omega$ is a ball centered at the origin. 

This circle of ideas has been developed in the last years. In \cite{AFS}, a Reilly-type inequality for the drift Laplacian was used to
derive a Heintze--Karcher inequality and a soap bubble theorem in weighted manifolds with
$\operatorname{Ric}^\alpha_w$ bounded from below, and the theorem of \cite{RWH_2018} was
extended to generalized cones. We
also mention the recent contributions \cite{AccCir,CFP,ADV,lima} on rigidity for overdetermined
and semilinear problems in weighted spaces. 


Throughout this paper, $d\mu_w=w\,d\mathrm{vol}$ and $d\sigma_w=w\,d\sigma$ are the
weighted volume and area elements, $H_w=H+\langle\nabla\log w,\nu\rangle$ is the weighted mean
curvature of $\partial\Omega$, and
\begin{equation*}
    r=\frac{Vol_w(\Omega)}{Area_w(\partial\Omega)},\qquad H_0=\frac1{r}.
\end{equation*}
In Theorems~\ref{main:A}, \ref{main:B} and~\ref{main:D}, $\Omega\subset\mathbb{R}^n$ is a
bounded domain with smooth boundary, $w$ is positive and homogeneous of degree $\alpha>0$, and
$q(x)=-\frac{|x|^2}{2(n+\alpha)}$. Our first result is the weighted counterpart of
\eqref{intro:MP}.

\begin{maintheorem}\label{main:A}
Let $u$ be the solution of \eqref{intro:problem}. Then, for every constant $c\in\mathbb{R}$,
\begin{equation}\label{intro:A}
    -2\int_\Omega u\Big[|\mathring{\nabla}^2u|^2
    +\operatorname{Ric}^\alpha_w(\nabla u,\nabla u)\Big]\,d\mu_w
    \geq\int_{\partial\Omega}\big(u_\nu^2-c^2\big)\big(u_\nu-q_\nu\big)\,d\sigma_w .
\end{equation}
In particular, if
\begin{equation*}
    \int_\Omega u\,\operatorname{Ric}^\alpha_w(\nabla u,\nabla u)\,d\mu_w\geq0
\end{equation*}
and the right-hand side of \eqref{intro:A} is nonnegative for some $c$, then $\Omega=B_R(0)$
and $u(x)=\frac{R^2-|x|^2}{2(n+\alpha)}$.
\end{maintheorem}

When $u_\nu$ is constant on $\partial\Omega$, the divergence theorem gives $u_\nu=-r$, and
the right-hand side of \eqref{intro:A} vanishes for $c=r$. Hence Theorem~\ref{main:A}
extends the weighted Serrin theorem of \cite{RWH_2018}: the curvature assumption is only
required in integral form, and the normalization $\nabla u(0)=0$ is not needed, since the fact
that the ball is centered at the origin is a consequence of the argument. Our second result is
a weighted version of the Magnanini--Poggesi formula for the Soap Bubble Theorem.

\begin{maintheorem}\label{main:B}
Let $u$ be the solution of \eqref{intro:problem} and assume that {{$0\notin\partial\Omega$}}. Then
\begin{align}\label{intro:B}
    \int_\Omega\Big[|\mathring{\nabla}^2u|^2+\operatorname{Ric}^\alpha_w(\nabla u,\nabla u)\Big]d\mu_w
    &+\frac{n+\alpha-1}{(n+\alpha)\,r}\int_{\partial\Omega}(u_\nu+r)^2\,d\sigma_w \nonumber\\
    &\leq\int_{\partial\Omega}\Big(\frac{n+\alpha-1}{n+\alpha}\,H_0-H_w\Big)(u_\nu-q_\nu)\,u_\nu\,d\sigma_w \nonumber\\
    &\quad+\int_{\partial\Omega}\Big(\frac{n+\alpha-1}{n+\alpha}\,H_0-H_w\Big)(u_\nu+r)\,q_\nu\,d\sigma_w .
\end{align}
In particular, if $\int_\Omega\operatorname{Ric}^\alpha_w(\nabla u,\nabla u)\,d\mu_w\geq0$ and
the right-hand side of \eqref{intro:B} is non-positive, then $\Omega$ is a ball centered at
the origin.
\end{maintheorem}

If $H_w$ is constant, the weighted Minkowski identity (Proposition~\ref{I-M}) gives
$H_w=\frac{n+\alpha-1}{n+\alpha}H_0$, the weighted mean curvature of the sphere of radius
$(n+\alpha)r$ centered at the origin. Then the right-hand side of \eqref{intro:B} vanishes,
and Theorem~\ref{main:B} yields a weighted Alexandrov theorem, which recovers
\cite[Corollary~2.9]{AFS} and shows, in addition, that the sphere is centered at the origin.
\begin{maincorollary}\label{main:cor}
Let $u$ be the solution of \eqref{intro:problem} and assume that
$\int_\Omega\operatorname{Ric}^\alpha_w(\nabla u,\nabla u)\,d\mu_w\geq0$. If $\partial\Omega$
has constant weighted mean curvature or, more generally, if
$H_w\geq\frac{n+\alpha-1}{n+\alpha}H_0$ on $\partial\Omega$, then $\Omega$ is a ball centered
at the origin.
\end{maincorollary}


Next, we turn to estimates for the weighted torsion function. Wang and Xia~\cite{wang2020}
obtained sharp bounds for the gradient and for the second normal derivative of the torsion
function of a compact Riemannian manifold with boundary whose Ricci curvature is bounded from
below; estimates of this type were recently obtained in \cite{AF}. These
estimates do not depend on the Euclidean structure nor on the homogeneity of the weight.

\begin{maintheorem}\label{main:C}
Let $\Omega$ be a bounded domain with smooth boundary in a weighted Riemannian manifold
$(M^n,g,w\,d\mathrm{vol})$, with $w$ smooth and positive on $\overline\Omega$, and let $u$ be
the solution of $\Delta_wu=-1$ in $\Omega$, $u=0$ on $\partial\Omega$. Assume that
$\operatorname{Ric}^\alpha_w\geq(n+\alpha-1)k\,g$ in $\Omega$, for some $\alpha>0$ and
$k\in\mathbb{R}$. Then
\begin{equation*}
    \max_{\partial\Omega}|\nabla u|^2>\frac{n+\alpha+2}{n+\alpha}\,\frac{T_w(\Omega)}{Vol_w(\Omega)}
    +\frac{2(n+\alpha-1)k}{Vol_w(\Omega)}\int_\Omega u\,|\nabla u|^2\,d\mu_w
\end{equation*}
and
\begin{equation*}
    \min_{\partial\Omega}u_{\nu\nu}<-\frac{1}{n+\alpha}-\frac{(n+\alpha-1)k\,T_w(\Omega)}{Vol_w(\Omega)} .
\end{equation*}
Moreover, if $H_w\geq0$ on $\partial\Omega$ and $H_w\not\equiv0$, then
\begin{equation*}
    \int_{\partial\Omega}u_{\nu\nu}\,d\sigma_w<(n+\alpha-1)^{1/2}
    \Big[\frac{Vol_w(\Omega)}{n+\alpha}-k\,T_w(\Omega)\Big]^{1/2}
    \Big[\int_{\partial\Omega}H_w\,d\sigma_w\Big]^{1/2}-Area_w(\partial\Omega).
\end{equation*}
\end{maintheorem}

In the unweighted case \cite{wang2020}, equality in these estimates holds if and only if $k=0$
and $\Omega$ is isometric to a Euclidean ball. Theorem~\ref{main:C} shows that this rigidity
disappears for smooth positive weights. 

Finally, we consider a \emph{reverse Serrin problem}. This problem was introduced by Magnanini,
Molinarolo and Poggesi~\cite{Magnanini2024}, who exchanged the roles of the two boundary
conditions in Serrin's problem: one starts from a solution of the \emph{Neumann} problem
\begin{equation}\label{intro:reverse}
    \Delta_wu=-1\quad\text{in }\Omega,\qquad u_\nu=-r\quad\text{on }\partial\Omega,
\end{equation}
which exists for every domain and is unique up to an additive constant, and asks which
conditions on the boundary values of $u$ force $\Omega$ to be a ball. In \cite{Magnanini2024}
this question is answered, in a quantitative form, by means of a general integral identity
involving two solutions of the Poisson equation, which contains \eqref{intro:MP} as a
particular case. We extend this identity to the drift Laplacian on arbitrary weighted
manifolds (Proposition~\ref{thm:general_weighted}) and deduce the following result, where
$u_*=\min_{\overline\Omega}u$, $h$ is the second fundamental form of $\partial\Omega$ and
$\bar\nabla$ is the tangential gradient on $\partial\Omega$.

\begin{maintheorem}\label{main:D}
Let $u\in C^3(\overline\Omega)$ be a solution of \eqref{intro:reverse} and let
$P=|\nabla u|^2+\frac{2u}{n+\alpha}$. Then
\begin{align}\label{intro:D}
    \int_\Omega(u-u_*)\,\Delta_wP\,d\mu_w
    &=\int_{\partial\Omega}|\bar\nabla u|^2(u_\nu-q_\nu)\,d\sigma_w
      -2\int_{\partial\Omega}(u-u_*)\,h(\bar\nabla u,\bar\nabla u)\,d\sigma_w \nonumber\\
    &\quad-2\int_{\partial\Omega}(u-u_*)\,(1-rH_w)(u_\nu-q_\nu)\,d\sigma_w .
\end{align}
In particular, if
$\int_\Omega(u-u_*)\operatorname{Ric}^\alpha_w(\nabla u,\nabla u)\,d\mu_w\geq0$ and the
right-hand side of \eqref{intro:D} is non-positive, then $\Omega$ is the ball of radius
$R=(n+\alpha)r$ centered at the origin, and
\begin{equation*}
    u(x)=u_*+\frac{R^2-|x|^2}{2(n+\alpha)} .
\end{equation*}
\end{maintheorem}

Every term on the right-hand side of \eqref{intro:D} is small when $u$ is close to a constant
on $\partial\Omega$; it vanishes when $u$ is constant on $\partial\Omega$, so that
Theorem~\ref{main:D} contains the weighted Serrin theorem (see
Corollary~\ref{cor:reverse_concrete}). 

\subsection{Organization of the paper}
In Section~2 we fix notation, recall the weighted Bochner and Reilly inequalities and the
weighted Minkowski identity, and collect the basic properties of homogeneous weights. Section~3
is devoted to the Dirichlet problem: we prove Theorems~\ref{main:A} and~\ref{main:B} (see
Theorems~\ref{thm:dirichlet}, \ref{thm:rigidity} and~\ref{thm:cone}, and
Corollary~\ref{cor:cone}) and the torsion estimates (Theorems~\ref{thm:max}). In Section~4 we establish the general weighted
identity and prove Theorem~\ref{main:D} (see Corollary~\ref{cor:weighted_3.3_u} and
Theorem~\ref{thm:reverse_serrin}).

\section{Preliminaries}


In this section we fix notation and recall the analytic tools for weighted
Riemannian manifolds that will be used throughout the paper. Here, $(M^n,\langle\cdot,\cdot\rangle, w\,d\mathrm{vol})$ denotes an
$n$-dimensional Riemannian manifold, $n\geq 2$, endowed with a smooth positive
weight $w\in C^\infty(M)$, and $d\mu_w = w\,d\mathrm{vol}$ denotes the
associated weighted measure. For a bounded domain $\Omega\subset M$ with smooth boundary
we write $Vol_w(\Omega)=\int_\Omega d\mu_w$ and
$Area_w(\partial\Omega)=\int_{\partial\Omega}d\sigma_w$, where $d\sigma_w=w\,d\sigma$. The
weighted divergence of a vector field $X$ is
$\operatorname{div}_wX=\operatorname{div}X+\langle\nabla\log w,X\rangle=w^{-1}\operatorname{div}(wX)$,
so that $\int_\Omega\operatorname{div}_wX\,d\mu_w=\int_{\partial\Omega}\langle X,\nu\rangle\,d\sigma_w$,
where $\nu$ is the outward unit normal. Finally, $\mathring{\nabla}^2u=\nabla^2u-\frac{\Delta u}{n}\,g$
denotes the traceless Hessian of $u$, so that
$|\nabla^2u|^2=|\mathring{\nabla}^2u|^2+\frac{(\Delta u)^2}{n}$.

The natural second-order operator associated with this structure is the
\emph{weighted Laplacian} (or $w$-Laplacian),
\begin{equation*}
    \Delta_w u = \Delta u + \langle \nabla \log w, \nabla u \rangle
    = \frac{1}{w}\operatorname{div}(w\nabla u).
\end{equation*}
Let $\Sigma^{n-1}$ be a hypersurface immersed in $M^n$ with unit normal $\nu$.
Intrinsic quantities on $\Sigma$ are denoted with an overline (e.g.\
$\bar\nabla$, $\bar{\Delta}_w$), where $\Sigma$ carries the
restricted weight $w|_{\Sigma}$ and the measure $d\sigma_w = w\,d\sigma$.
We denote by $h(X,Y)=\langle \nabla_X\nu, Y\rangle$ the second fundamental
form of $\Sigma$, by $H=\operatorname{tr} h=\operatorname{div}_\Sigma \nu$ its
mean curvature, and by
\begin{equation*}
    H_w = H + \langle \nabla \log w, \nu \rangle
\end{equation*}
its weighted (Gromov) mean curvature \cite{gromov}.

The following lemma relates the ambient and intrinsic weighted Laplacians;
see, for instance, \cite[Lemma 2.1]{AFS}.

\begin{lemma}\label{lem:laplacians}
    Let $\Sigma^{n-1}$ be a hypersurface immersed in $M^n$ and let
    $f\in C^\infty(M)$. Then, along $\Sigma$,
    \begin{equation*}
        \Delta_w f = \bar{\Delta}_w f + H_w f_\nu + \nabla^2 f(\nu,\nu),
    \end{equation*}
    where $f_\nu = \langle \nabla f, \nu\rangle$ and $\nabla^2 f$ denotes the
    ambient Hessian of $f$.
\end{lemma}

Next we recall the weighted Bochner inequality, which plays a central role in
our work. It follows by combining the classical Bochner
formula with the Cauchy--Schwarz inequality and Bergström's inequality
\cite{ber}; see \cite[Appendix~A]{wey}, and \cite[Section~1.5]{li} for an
elementary proof.

\begin{lemma}\label{Weightedbochner}
    Let $(M^n,\langle\cdot,\cdot\rangle,w\,d\mathrm{vol})$ be a weighted
    Riemannian manifold, let $u\in C^\infty(M)$ and let $\alpha>0$. Then
    \begin{equation*}
        \frac{1}{2}\Delta_w |\nabla u|^2
        \geq \frac{(\Delta_w u)^2}{n+\alpha}
        + \langle \nabla u,\nabla\Delta_w u\rangle
        + \operatorname{Ric}_{w}^{\alpha}(\nabla u,\nabla u),
    \end{equation*}
    where
    \begin{equation*}
        \operatorname{Ric}_{w}^{\alpha}
        = \operatorname{Ric} - \nabla^2 \log w
        - \frac{1}{\alpha}\, d\log w \otimes d\log w
    \end{equation*}
    is the $\alpha$-Bakry--Émery--Ricci tensor. Moreover, equality holds at a
    point if and only if, at that point,
    \begin{equation}
        \nabla^2 u = \frac{\Delta u}{n}\, g
        \quad\text{and}\quad
        \frac{\Delta u}{n}
        = \frac{\langle\nabla\log w,\nabla u\rangle}{\alpha}
        = \frac{\Delta_w u}{n+\alpha}\label{nablaunablaw}.
    \end{equation}
\end{lemma}


Integrating Lemma~\ref{Weightedbochner} over a compact manifold with boundary
yields the following weighted Reilly-type inequality, together with a
characterization of the equality case. It was first established in
\cite[Proposition~2.3]{AFS}.

\begin{proposition}\label{Reilly}
    Let $(M^n,\langle\cdot,\cdot\rangle,w\,d\mathrm{vol})$ be a compact
    weighted Riemannian manifold with smooth boundary $\partial M$, let
    $u\in C^\infty(\overline{M})$ and let $\alpha>0$. Set $z=u|_{\partial M}$
    and let $\nu$ be the outward unit normal to $\partial M$. Then
    \begin{align}
        \frac{n+\alpha-1}{n+\alpha}\int_{M}(\Delta_w u)^2\,d\mu_w
        &\geq \int_{M}\operatorname{Ric}_{w}^{\alpha}(\nabla u,\nabla u)\,d\mu_w
        \nonumber\\
        &\quad + \int_{\partial M}
        \Big[\,2u_{\nu}\,\bar{\Delta}_w z + H_w u_{\nu}^2
        + h(\bar\nabla z,\bar\nabla z)\Big]\,d\sigma_w,
        \label{eqreilly}
    \end{align}
    where $H_w$ and $h$ are the weighted mean curvature and the second
    fundamental form of $\partial M$ with respect to $\nu$. Moreover, equality
    holds in \eqref{eqreilly} if and only if
    \begin{equation*}
        \nabla^2 u = \frac{\Delta u}{n}\, g
        \quad\text{and}\quad
        \frac{\Delta u}{n}
        = \frac{\langle\nabla\log w,\nabla u\rangle}{\alpha}
        = \frac{\Delta_w u}{n+\alpha}
        \quad\text{in } M.
    \end{equation*}
\end{proposition}

Finally, we recall a weighted version of the Minkowski identity for weights
that are positive and homogeneous of degree $\alpha>0$, which means
$\langle\nabla w,x\rangle=\alpha w$, where $x$ is the position vector. Such a weight is
defined in $\mathbb{R}^n\setminus\{0\}$ and satisfies $w(x)=|x|^\alpha\,w(x/|x|)\to0$ as
$x\to0$, so the origin is a singularity of the weighted structure; see
\cite[Remark~3.2]{CC_2014}. For this reason, hypersurfaces are always assumed to lie in
$\mathbb{R}^n\setminus\{0\}$, whereas domains are allowed to contain the origin. Versions
of this identity for hypersurfaces with nonempty boundary can be found in
\cite[Proposition~5.1]{CC_2014}.

\begin{proposition}\label{I-M}
    Let $w\in C^1(\mathbb{R}^n\setminus\{0\})$ be a positive function,
    homogeneous of degree $\alpha>0$, and let
    $\Sigma^{n-1}\subset\mathbb{R}^n\setminus\{0\}$ be a closed, orientable
    hypersurface with unit normal $\nu$. Then
    \begin{equation*}
        \int_{\Sigma}\Big[(n+\alpha-1) - H_w\langle x,\nu\rangle\Big]\,d\sigma_w = 0.
    \end{equation*}
\end{proposition}

We collect the properties of homogeneous weights that will be used repeatedly.

\begin{lemma}\label{lem:homogeneous}
Let $w\in C^2(\mathbb{R}^n\setminus\{0\})$ be positive and homogeneous of degree $\alpha>0$,
and set $f=w^{1/\alpha}$. Then the following hold in $\mathbb{R}^n\setminus\{0\}$.
\begin{enumerate}[(a)]
\item $\langle\nabla\log w,x\rangle=\alpha$ and $\operatorname{div}_wx=n+\alpha$. Consequently,
the function $q(x)=-\frac{|x|^2}{2(n+\alpha)}$ satisfies
\begin{equation*}
    \Delta_wq=-1,\qquad\nabla^2q=-\frac{1}{n+\alpha}\,g
    \qquad\text{and}\qquad
    \frac{\Delta q}{n}=\frac{\langle\nabla\log w,\nabla q\rangle}{\alpha}=-\frac{1}{n+\alpha}.
\end{equation*}
In particular, equality holds in Lemma~\ref{Weightedbochner} for $u=q$. Moreover,
$q_\nu=-\frac{\langle x,\nu\rangle}{n+\alpha}$ along any hypersurface with unit normal $\nu$.
\item $\operatorname{Ric}^\alpha_w=-\dfrac{\alpha}{f}\,\nabla^2f$.
\item $\operatorname{Ric}^\alpha_w(x,\cdot)=0$.
\end{enumerate}
\end{lemma}

\begin{proof}
Item (a) follows from a direct computation. For (b) and (c), recall that $\log w=\alpha\log f$.
Hence
\begin{equation*}
    \nabla^2\log w=\alpha\Big(\frac{\nabla^2f}{f}-\frac{df\otimes df}{f^2}\Big)
    \qquad\text{and}\qquad
    \frac1\alpha\,d\log w\otimes d\log w=\alpha\,\frac{df\otimes df}{f^2},
\end{equation*}
and (b) follows because the Ricci tensor of $\mathbb{R}^n$ vanishes. Moreover, Euler's
identity $\langle\nabla w,x\rangle=\alpha w$ becomes $\langle\nabla f,x\rangle=f$.
Differentiating it in the direction of a vector $Y$ and using $\nabla_Yx=Y$, we get
$\nabla^2f(x,Y)=0$, and (c) follows from (b).
\end{proof}

\begin{remark}\label{rem:origin}
In our Euclidean results, $w\in C^2(\mathbb{R}^n\setminus\{0\})$ is positive and homogeneous
of degree $\alpha>0$, and the domain $\Omega$ is allowed to contain the origin. Since $w$ is
bounded near the origin (see \cite[Remark~4.1]{CC_2014}), the integrations by parts used
below remain valid in this case, as observed in \cite[Remark~1.1]{RWH_2018}.
\end{remark}

\section{Torsion Problem for the Drift Laplacian and rigidity}

\subsection{Spherical deficit and rigidity}

In this subsection, $\Omega\subset\mathbb{R}^n$ and $w$ satisfy the standing assumptions of
Section~2. We consider the following boundary value problem:
\begin{equation}\label{problem1}
\left\{\begin{array}{rcl}
\Delta_\omega u  &=& -1\quad \text{in}\quad  \Omega,\\
u&=& 0\quad \text{on}\quad  \partial\Omega.
\end{array}\right.
\end{equation}
By the maximum principle and the Hopf lemma \cite{GT}, $u>0$ in $\Omega$ and $u_\nu<0$ on
$\partial\Omega$.

Our goal here is to study the geometric properties of the solutions to
\eqref{problem1}. To this end, we introduce an auxiliary function $P$ (a
generalization of the classical $P$-function of Weinberger~\cite{Weinberger}), defined by
\begin{equation}\label{PFunction}
    P=|\nabla u|^2+\frac{2u}{n+\alpha}.
\end{equation}
This choice is motivated by its usefulness in applying maximum principles and
integral identities. By applying the weighted Bochner formula
(Lemma~\ref{Weightedbochner}) to this $P$-function, we obtain a differential
inequality that governs its behavior. This is the content of our next lemma,
which will be fundamental for establishing our rigidity results and estimates
for the spherical deficit in the weighted setting (see \cite[Lemma~2.2]{AFS}).
The lemma holds on any weighted manifold.

\begin{lemma}\label{BochnerP}
Let $(M^n, \langle\cdot,\cdot\rangle, w\,d\mathrm{vol})$ be an $n$-dimensional
weighted manifold, $n\geq 2$. Let $\Omega\subset M$ be a domain and let
$u\in C^3(\Omega)$ be a solution to $\Delta_w u=-1$. Then
\[
\Delta_{w}P\geq 2\left[|\mathring{\nabla}^2 u|^2+\operatorname{Ric}_w^{\alpha}(\nabla u,\nabla u)\right].
\]
Equality holds if and only if
$\frac{\Delta u}{n} = \frac{\langle \nabla \log w, \nabla u \rangle}{\alpha} = \frac{\Delta_w u}{n+\alpha}$.
\end{lemma}


The next lemma describes the solutions for which both the spherical deficit vanishes and
equality holds in Lemma~\ref{BochnerP}. It will be used in all our rigidity statements.

\begin{lemma}\label{lem:quadratic}
Let $\Omega\subset\mathbb{R}^n$ and $w$ satisfy the standing assumptions, and let $u$ be a
solution of $\Delta_wu=-1$ such that, in $\Omega\setminus\{0\}$,
\begin{equation*}
    \mathring{\nabla}^2u=0\qquad\text{and}\qquad
    \frac{\Delta u}{n}=\frac{\langle\nabla\log w,\nabla u\rangle}{\alpha}.
\end{equation*}
Then there exist $A\in\mathbb{R}$ and $z\in\mathbb{R}^n$ such that
\begin{equation}\label{eq:quadratic}
    u(x)=A-\frac{|x-z|^2}{2(n+\alpha)}\quad\text{in }\Omega
    \qquad\text{and}\qquad
    \langle\nabla\log w,x-z\rangle=\alpha\quad\text{in }\Omega\setminus\{0\}.
\end{equation}
Moreover, $\operatorname{Ric}^\alpha_w(\nabla u,\nabla u)=0$ in $\Omega\setminus\{0\}$, and if
$z\in\Omega$, then $z=0$. In particular, if $u>0$ in $\Omega$ and $u=0$ on $\partial\Omega$,
then $\Omega=B_R(0)$ and $u(x)=\frac{R^2-|x|^2}{2(n+\alpha)}$.
\end{lemma}

\begin{proof}
By assumption, $\frac{\Delta u}{n}=\frac{\langle\nabla\log w,\nabla u\rangle}{\alpha}
=\frac{\Delta_wu}{n+\alpha}=-\frac{1}{n+\alpha}$. Together with $\mathring{\nabla}^2u=0$, this
gives $\nabla^2u=-\frac{1}{n+\alpha}\,g$ in $\Omega\setminus\{0\}$, which is connected because
$n\geq2$. Hence $u$ has the form \eqref{eq:quadratic}, and $\nabla u=-\frac{x-z}{n+\alpha}$.
The identity $\langle\nabla\log w,\nabla u\rangle=-\frac{\alpha}{n+\alpha}$ then yields
$\langle\nabla\log w,x-z\rangle=\alpha$.

Comparing with $\langle\nabla\log w,x\rangle=\alpha$, we get $\langle\nabla f,z\rangle=0$ in
$\Omega\setminus\{0\}$, where $f=w^{1/\alpha}$. Differentiating, $\nabla^2f(z,\cdot)=0$, and
since $\nabla^2f(x,\cdot)=0$ we obtain $\nabla^2f(x-z,\cdot)=0$. By
Lemma~\ref{lem:homogeneous}(b), $\operatorname{Ric}^\alpha_w(\nabla u,\nabla u)=0$.

If $z\in\Omega$ and $z\neq0$, then $w$ is of class $C^1$ and positive near $z$, and evaluating
$\langle\nabla\log w,x-z\rangle=\alpha$ at $x=z$ gives $0=\alpha$, a contradiction. Finally,
if $u>0$ in $\Omega$ and $u=0$ on $\partial\Omega$, then $A>0$ and $\Omega$ is the ball of radius
$R=\sqrt{2(n+\alpha)A}$ centered at $z$; in particular $z\in\Omega$, so $z=0$.
\end{proof}

We also see a weighted Pohozaev identity, which will be used in the proof of
Theorem~\ref{thm:dirichlet}; see also \cite{pohozaev}, and \cite[Proposition~3.2]{AFS}.

\begin{lemma}\label{lem:pohozaev}
Let $u$ be the solution of \eqref{problem1}. Then
\begin{equation*}
    (n+\alpha+2)\int_\Omega|\nabla u|^2\,d\mu_w=\int_{\partial\Omega}\langle x,\nu\rangle\,u_\nu^2\,d\sigma_w .
\end{equation*}
\end{lemma}

\begin{proof}
Let $W = \langle x, \nabla u \rangle \nabla u$. Using the
definition of the weighted divergence
\[
\operatorname{div}_\omega W = \operatorname{div} W + \langle \nabla \log \omega, W \rangle,
\]
we obtain
\begin{equation*}
    \operatorname{div}_\omega W = - \langle x, \nabla u \rangle + |\nabla u|^2 + \frac{1}{2} \langle x, \nabla |\nabla u|^2 \rangle,
\end{equation*}
where we used that $\Delta_\omega u = -1$. Then, from the divergence theorem,
and since $\nabla u=u_\nu\nu$ on $\partial\Omega$, it follows that
\begin{equation}\label{div}
   \int_{\partial \Omega} \langle x, \nabla u \rangle \langle \nabla u, \nu \rangle \,d\sigma_\omega
   = \int_{\Omega} \left( - \langle x, \nabla u \rangle + |\nabla u|^2 + \frac{1}{2} \langle x, \nabla |\nabla u|^2 \rangle \right) \,d\mu_\omega.
 \end{equation}

Now, using that $\operatorname{div}_\omega(ux)=\langle\nabla u,x\rangle+u\operatorname{div}_\omega(x)$
and $\operatorname{div}_\omega(x|\nabla u|^2)=(n+\alpha)|\nabla u|^2+\langle x,\nabla|\nabla u|^2\rangle$,
with $\operatorname{div}_\omega(x)=n+\alpha$, from the divergence theorem we get
\[
\int_\Omega\langle\nabla u,x\rangle\,d\mu_\omega = -\int_\Omega u(n+\alpha)\,d\mu_\omega
\]
and
\[
\int_\Omega\langle x,\nabla|\nabla u|^2\rangle\,d\mu_\omega = -\int_\Omega(n+\alpha)|\nabla u|^2\,d\mu_\omega+\int_{\partial\Omega}|\nabla u|^2\langle x,\nu\rangle\,d\sigma_\omega.
\]
Thus, \eqref{div} becomes
\begin{align*}
\int_{\partial \Omega} \langle x, \nu \rangle \, u_\nu^2 \,d\sigma_\omega
= \int_{\Omega} \left( u(n+\alpha) + |\nabla u|^2 - \frac{1}{2} (n+\alpha) |\nabla u|^2 \right) \,d\mu_\omega
+ \frac{1}{2} \int_{\partial \Omega} u_\nu^2 \langle x, \nu \rangle \,d\sigma_\omega.
\end{align*}
Therefore,
\begin{align*}
    \frac{1}{2} \int_{\partial \Omega} \langle x, \nu \rangle \, u_\nu^2 \,d\sigma_\omega
    &= (n+\alpha) \int_{\Omega} u \,d\mu_\omega + \frac{2 - (n+\alpha)}{2} \int_{\Omega} |\nabla u|^2 \,d\mu_\omega\\
    &= (n+\alpha) \int_{\Omega} |\nabla u|^2 \,d\mu_\omega + \frac{2 - (n+\alpha)}{2} \int_{\Omega} |\nabla u|^2 \,d\mu_\omega\\
    &= \frac{(n+\alpha+2)}{2} \int_{\Omega} |\nabla u|^2 \,d\mu_\omega,
\end{align*}
where we have used that
$\int_{\Omega} u\,d\mu_{\omega}=\int_{\Omega}|\nabla u|^2\,d\mu_{\omega}.$
\end{proof}

We are now in a position to state the next result, which was inspired by
\cite[Theorem~2.1]{Magnanini2020}. By integrating the inequality provided by
Lemma~\ref{BochnerP} and comparing the solution $u$ with a suitable quadratic
comparison function $q$, we can establish a sharp integral inequality relating
the interior geometry of $\Omega$ to its boundary data.

\begin{theorem}\label{thm:dirichlet}
    Let $\Omega\subset\mathbb{R}^n$ be a bounded domain, let $w$ be positive and
    homogeneous of degree $\alpha>0$, and let $c\in\mathbb{R}$ be a constant.
    Then, the solution $u$ of \eqref{problem1} satisfies
\begin{equation}\label{Dirichlet_inequality}
    -\int_\Omega u\left[2|\mathring{\nabla}^2 u|^2+2\operatorname{Ric}_\omega^{\alpha}(\nabla u,\nabla u)\right] \, d\mu_\omega
    \geq\int_{\partial\Omega}(u_\nu^2-c^2)(u_\nu -q_\nu) \, d\sigma_\omega,
\end{equation}
    where $\nu$ denotes the outward unit normal vector on $\partial\Omega$, and
    $q$ is the comparison solution of
    \[
    \Delta_\omega q=-1,
    \]
    given by
    \[
    q(x)=-\frac{|x|^2}{2(n+\alpha)}.
    \]
\end{theorem}

\begin{proof}
The proof is divided into two main steps. First, we use the divergence theorem
and the properties of the homogeneous weight to express the boundary terms in
terms of the comparison function $q$. Subsequently, we apply the Bochner
inequality to bound the volume integral.

We begin by using the fact that $u$ is a solution of \eqref{problem1}, which
yields
\begin{equation*}
    \int_\Omega (u\Delta_\omega P-P\Delta_\omega u)\,d\mu_\omega
    =\int_{\partial\Omega}(uP_\nu-Pu_\nu)\,d\sigma_\omega
    = -\int_{\partial\Omega}u_\nu^3\,d\sigma_\omega,
\end{equation*}
where in the last equality we used that $u=0$ and $P=u_\nu^2$ on $\partial\Omega$.
Then,
\begin{equation}\label{laplacian P}
    -\int_\Omega u\Delta_\omega P\,d\mu_\omega
    =\int_\Omega P\,d\mu_\omega+\int_{\partial\Omega}u_\nu^3\,d\sigma_\omega.
\end{equation}

Now, observe that
\begin{align}\label{int}
    \int_\Omega P\,d\mu_\omega
    &=\int_\Omega |\nabla u|^2\,d\mu_\omega-\frac{2}{n+\alpha}\int_\Omega u\Delta_\omega u\,d\mu_\omega \nonumber\\
    &=\left[1+\frac{2}{n+\alpha}\right]\int_\Omega |\nabla u|^2\,d\mu_\omega.
\end{align}
Back to \eqref{int}, by Lemma~\ref{lem:pohozaev} we get
\[
\int_\Omega P\,d\mu_\omega = \frac{1}{n+\alpha}\int_{\partial\Omega} \langle x,\nu\rangle u_\nu^2\,d\sigma_\omega.
\]
Then, from \eqref{laplacian P}, and since $q_\nu=-\frac{\langle x,\nu\rangle}{n+\alpha}$
(Lemma~\ref{lem:homogeneous}(a)), we have
\begin{align*}
      -\int_\Omega u\Delta_\omega P\,d\mu_\omega
      &=\frac{1}{n+\alpha}\int_{\partial\Omega} \langle x,\nu\rangle u_\nu^2\,d\sigma_\omega+\int_{\partial\Omega}u_\nu^3\,d\sigma_\omega\\
      &=-\int_{\partial\Omega} q_\nu u_\nu^2\,d\sigma_\omega+\int_{\partial\Omega}u_\nu^3\,d\sigma_\omega\\
      &=-\int_{\partial\Omega} q_\nu u_\nu^2\,d\sigma_\omega+\int_{\partial\Omega}\left[(u_\nu^2-c^2)(u_\nu-q_\nu)+c^2(u_\nu-q_\nu)
      + u_\nu^2q_\nu\right]\,d\sigma_\omega\\
      &=\int_{\partial\Omega}(u_\nu^2-c^2)(u_\nu-q_\nu)\,d\sigma_\omega,
\end{align*}
where we used that
\[
0=\int_\Omega\Delta_\omega(u-q)\,d\mu_\omega=\int_{\partial\Omega}(u_\nu-q_\nu)\,d\sigma_\omega.
\]

On the other hand, from Lemma~\ref{BochnerP} we have that
\[
-\int_\Omega u\Delta_\omega P\,d\mu_\omega
\leq-\int_\Omega u\left[2|\mathring{\nabla}^2 u|^2+2\operatorname{Ric}_\omega^{\alpha}(\nabla u,\nabla u)\right] \,d\mu_\omega,
\]
where we have used that $u>0$ (by the maximum principle). Finally, we obtain
\begin{equation*}
-2\int_\Omega u\left[|\mathring{\nabla}^2 u|^2+\operatorname{Ric}_\omega^{\alpha}(\nabla u,\nabla u)\right] \,d\mu_\omega
\geq\int_{\partial\Omega}(u_\nu^2-c^2)(u_\nu -q_\nu) \,d\sigma_\omega.
\end{equation*}
\end{proof}

\begin{theorem}\label{thm:rigidity}
Let $\Omega$ and $w$ be as in Theorem~\ref{thm:dirichlet}, and let $u$ be the solution of
\eqref{problem1}. If the right-hand side of \eqref{Dirichlet_inequality} is nonnegative for
some $c\in\mathbb{R}$ and

\begin{equation}\label{eq:int_ric_u}
    \int_{\Omega}u\,\operatorname{Ric}_\omega^{\alpha}(\nabla u,\nabla u)\,d\mu_\omega\geq 0,
\end{equation}
then $\Omega=B_R(0)$ for some $R>0$, and $u(x)=\frac{R^2-|x|^2}{2(n+\alpha)}$.
\end{theorem}


\begin{proof}
By the hypotheses, the left-hand side of \eqref{Dirichlet_inequality} is at most
$-2\int_\Omega u\,|\mathring{\nabla}^2u|^2\,d\mu_w\leq0$, while its right-hand side is
nonnegative. Hence $\int_\Omega u\,|\mathring{\nabla}^2u|^2\,d\mu_w=0$, and all the inequalities
in the proof of Theorem~\ref{thm:dirichlet} are equalities. In particular,
\begin{equation*}
    \int_\Omega u\Big(\Delta_wP-2\big[|\mathring{\nabla}^2u|^2+\operatorname{Ric}^\alpha_w(\nabla u,\nabla u)\big]\Big)\,d\mu_w=0,
\end{equation*}
where the integrand is nonnegative by Lemma~\ref{BochnerP}. Since $u>0$ in $\Omega$, we get
$\mathring{\nabla}^2u=0$ and equality in Lemma~\ref{BochnerP} at every point of
$\Omega\setminus\{0\}$. The conclusion follows from Lemma~\ref{lem:quadratic}.
\end{proof}


\begin{remark}
The previous result gives an alternative proof for \cite[Theorem~1.1]{RWH_2018}.
We point out that it is not necessary to assume the condition $\nabla u(0)=0$. Indeed, if $u_\nu$ is constant on $\partial\Omega$, then $u_\nu=-r$ with
$r=Vol_w(\Omega)/Area_w(\partial\Omega)$ by the divergence theorem, and the right-hand side
of \eqref{Dirichlet_inequality} vanishes for $c=r$. Moreover, the curvature assumption is
only required in the integral form \eqref{eq:int_ric_u}.
\end{remark}


Now, we derive an integral identity relating the interior geometry of $\Omega$
to the boundary curvatures to get the following result.

\begin{theorem}\label{thm:cone}
Let $\Omega\subset\mathbb{R}^n$ and $w$ be as in Theorem~\ref{thm:dirichlet}, with $0\notin\partial\Omega$. If $u$ is a solution to \eqref{problem1} and $\Omega$ is a domain whose
boundary $\partial\Omega$ satisfies the conditions of Proposition~\ref{I-M},
then
\begin{align}\label{thm-M}
&\int_{\Omega}[|\mathring{\nabla}^2 u|^2+\operatorname{Ric}^{\alpha}_w(\nabla u, \nabla u)]\,d\mu_\omega+\frac{(n+\alpha-1)}{(n+\alpha)r}\int_{\partial \Omega}(u_{\nu}+r)^2\,d\sigma_\omega\nonumber\\
&\leq\int_{\partial \Omega}\left(\frac{n+\alpha-1}{n+\alpha}H_0-H_w\right)(u_{\nu}-q_{\nu})u_{\nu}\,d\sigma_\omega\nonumber\\
&\quad+\int_{\partial \Omega}\left(\frac{n+\alpha-1}{n+\alpha}H_0-H_w\right)(u_{\nu}+r)q_{\nu}\,d\sigma_\omega,
\end{align}
where $r=\frac{Vol_w(\Omega)}{Area_w(\partial \Omega)}$, $H_0=1/r$ and
$q(x)=-\frac{|x|^2}{2(n+\alpha)}$.
\end{theorem}


\begin{proof}
Since $P=|\nabla u|^2+\frac{2}{n+\alpha}u$, we compute its normal derivative:
\[
P_{\nu}=\langle\nabla|\nabla u|^2,\nu\rangle+\frac{2}{n+\alpha}u_{\nu}=2\nabla^2u(\nabla u, \nu)+\frac{2}{n+\alpha}u_{\nu}.
\]
Integrating by parts, we obtain
\[
\int_{\Omega}\Delta_wP\,d\mu_\omega=\int_{\partial\Omega}P_{\nu}\,d\sigma_\omega.
\]
Next, using the Bochner inequality in Lemma~\ref{BochnerP}, and since $\nabla u=u_\nu\nu$
on $\partial\Omega$, we obtain
\begin{equation}\label{eq006}
\begin{aligned}
\int_{\Omega}\left(|\mathring{\nabla}^2 u|^2+\operatorname{Ric}^{\alpha}_w(\nabla u,\nabla u)\right)d\mu_\omega
&\leq\int_{\partial\Omega}\left[\nabla^2u(\nabla u,\nu)+\frac{1}{n+\alpha}u_{\nu}\right]d\sigma_\omega\\
&=\int_{\partial\Omega}\left[u_\nu\nabla^2u(\nu,\nu)+\frac{1}{n+\alpha}u_{\nu}\right]d\sigma_\omega.
\end{aligned}
\end{equation}
At the same time, using that $\Delta_wu=\bar{\Delta}_wu+H_wu_{\nu}+\nabla^2u(\nu,\nu)$, we obtain
$\nabla^2u(\nu,\nu)=-1-H_wu_{\nu}$ on $\partial\Omega$. Combining this with \eqref{eq006}, we infer
\begin{equation}\label{eq007}
\begin{aligned}
\int_{\partial\Omega}\left[u_\nu\nabla^2u(\nu,\nu)+\frac{1}{n+\alpha}u_{\nu}\right]d\sigma_\omega
&=\int_{\partial \Omega}u_{\nu}\left[-1-H_wu_{\nu}+\frac{1}{n+\alpha}\right]d\sigma_\omega\\
&=\int_{\partial \Omega}\left[-H_wu_{\nu}^2+\frac{(-n-\alpha+1)}{n+\alpha}u_{\nu}\right]d\sigma_\omega.
\end{aligned}
\end{equation}
Now, since $H_0=\frac{1}{r}$, we have
\begin{equation}\label{eq008}
\begin{aligned}
\frac{1}{r}\int_{\partial \Omega}(u_{\nu}+r)^2\,d\sigma_\omega
&=\frac{1}{r}\int_{\partial \Omega}(u^2_{\nu}+2u_{\nu}r+r^2)\,d\sigma_\omega\\
&=H_0\int_{\partial \Omega}u^2_{\nu}\,d\sigma_\omega+2\int_{\partial \Omega}u_{\nu}\,d\sigma_\omega+r\,Area_w(\partial\Omega)\\
&=H_0\int_{\partial \Omega}u^2_{\nu}\,d\sigma_\omega+\int_{\partial \Omega}u_{\nu}\,d\sigma_\omega,
\end{aligned}
\end{equation}
because $r Area_w(\partial\Omega)=-\int_{\partial\Omega}u_{\nu}\,d\sigma_\omega=Vol_w(\Omega)$.
This implies that
\begin{equation}\label{eq009}
\frac{(n+\alpha-1)}{r(n+\alpha)}\int_{\partial \Omega}(u_{\nu}+r)^2\,d\sigma_\omega
=\frac{(n+\alpha-1)}{n+\alpha}\left(H_0\int_{\partial \Omega}u^2_{\nu}\,d\sigma_\omega+\int_{\partial \Omega}u_{\nu}\,d\sigma_\omega\right).
\end{equation}
From \eqref{eq006}, \eqref{eq007} and \eqref{eq009}, we conclude
\begin{equation}\label{eq:soap_intermediate}
\begin{aligned}
&\int_{\Omega}[|\mathring{\nabla}^2 u|^2+\operatorname{Ric}^{\alpha}_w(\nabla u,\nabla u)]\,d\mu_\omega+\frac{n+\alpha -1}{(n+\alpha)r}\int_{\partial\Omega}(u_{\nu}+r)^2\,d\sigma_\omega\\
&\leq \int_{\partial\Omega}u^2_{\nu}\left[-H_w+\frac{(n+\alpha-1)}{n+\alpha}H_0\right]d\sigma_\omega.
\end{aligned}
\end{equation}
Now, writing $u_\nu^2=(u_\nu-q_\nu)u_\nu+(u_\nu+r)q_\nu-r q_\nu$, we get
\begin{align*}
&\int_{\Omega}[|\mathring{\nabla}^2 u|^2+\operatorname{Ric}^{\alpha}_w(\nabla u,\nabla u)]\,d\mu_\omega+\frac{n+\alpha -1}{(n+\alpha)r}\int_{\partial\Omega}(u_{\nu}+r)^2\,d\sigma_\omega\\
&\leq \int_{\partial\Omega}\left[-H_w+\frac{(n+\alpha-1)}{n+\alpha}H_0\right](u_\nu-q_\nu)u_\nu\,d\sigma_\omega\\
&\quad+ \int_{\partial\Omega}\left[-H_w+\frac{(n+\alpha-1)}{n+\alpha}H_0\right](u_\nu+r)q_\nu\,d\sigma_\omega
- \int_{\partial\Omega}\left[-H_w+\frac{(n+\alpha-1)}{n+\alpha}H_0\right]r q_\nu\,d\sigma_\omega.
\end{align*}
Moreover,
\[
\int_{\partial\Omega}\left[\left(\frac{n+\alpha-1}{n+\alpha}\right)H_0-H_w\right]r q_{\nu}\,d\sigma_\omega=0.
\]
Indeed, since $q(x)=-\frac{|x|^2}{2(n+\alpha)}$, we have $\Delta_wq=-1$ and
$q_{\nu}=-\frac{\langle x,\nu\rangle}{n+\alpha}$. Integrating by parts, we obtain
\[
\int_{\partial\Omega}q_{\nu}\,d\sigma_\omega=\int_{\Omega}\Delta_wq\,d\mu_\omega=-Vol_w(\Omega).
\]
Now, using Proposition~\ref{I-M}, we get
\[
Area_w(\partial\Omega)\frac{[n+\alpha-1]}{n+\alpha}+\int_{\partial\Omega}H_wq_{\nu}\,d\sigma_\omega=0.
\]
Therefore,
\[
r\int_{\partial \Omega}H_wq_{\nu}\,d\sigma_\omega
=-\frac{r(n+\alpha-1)}{n+\alpha}Area_w(\partial \Omega)
=\frac{(n+\alpha-1)}{n+\alpha}\int_{\partial\Omega}q_{\nu}\,d\sigma_\omega,
\]
and the claim follows because $H_0r=1$.
This completes the proof.
\end{proof}

\begin{corollary}\label{cor:cone}
Suppose $u$ is a solution of \eqref{problem1}. If the right-hand side of
\eqref{thm-M} is non-positive and
$\displaystyle\int_{\Omega}\operatorname{Ric}_\omega^{\alpha}(\nabla u,\nabla u)\,d\mu_\omega\geq 0$, then $\Omega$
is a ball centered at the origin. In particular, this is the case if
$H_w\geq\frac{n+\alpha-1}{n+\alpha}H_0$ on $\partial\Omega$, for instance if $H_w$ is
constant on $\partial\Omega$.
\end{corollary}


\begin{proof}
Under the assumptions, the left-hand side of \eqref{thm-M} is at least
$\int_\Omega|\mathring{\nabla}^2u|^2\,d\mu_w\geq0$, and it is bounded above by a non-positive
quantity. Hence $\mathring{\nabla}^2u=0$, and all the inequalities in the proof of
Theorem~\ref{thm:cone} are equalities. In particular,
$\int_\Omega\big(\Delta_wP-2\big[|\mathring{\nabla}^2u|^2+\operatorname{Ric}^\alpha_w(\nabla
u,\nabla u)\big]\big)\,d\mu_w=0$ with a nonnegative integrand, so that equality holds in
Lemma~\ref{BochnerP} at every point of $\Omega\setminus\{0\}$. The conclusion follows from
Lemma~\ref{lem:quadratic}.

For the last assertion, recall from \eqref{eq:soap_intermediate} and the end of the proof of
Theorem~\ref{thm:cone} that the right-hand side of \eqref{thm-M} equals
$\int_{\partial\Omega}\big(\frac{n+\alpha-1}{n+\alpha}H_0-H_w\big)u_\nu^2\,d\sigma_w$, which
is non-positive if $H_w\geq\frac{n+\alpha-1}{n+\alpha}H_0$. Finally, if $H_w$ is constant,
Proposition~\ref{I-M} and the divergence theorem give
\[
(n+\alpha-1)\,Area_w(\partial\Omega)=H_w\int_{\partial\Omega}\langle x,\nu\rangle\,d\sigma_w
=H_w\int_\Omega\operatorname{div}_wx\,d\mu_w=(n+\alpha)\,H_w\,Vol_w(\Omega),
\]
that is, $H_w=\frac{n+\alpha-1}{n+\alpha}H_0$.
\end{proof}

\subsection{Estimates of the weighted torsion function}
The estimates of this subsection do not rely on the Euclidean structure nor on the
homogeneity of the weight, so we first state them on a general weighted manifold. Let $\Omega$
be a bounded domain with smooth boundary in a weighted Riemannian manifold
$(M^n,\langle\cdot,\cdot\rangle,w\,d\mathrm{vol})$, with $w$ smooth and positive on
$\overline\Omega$, and let $u\in C^3(\overline\Omega)$ be the solution of
\begin{equation}\label{eq:torsion_general}
    \Delta_wu=-1\quad\text{in }\Omega,\qquad u=0\quad\text{on }\partial\Omega .
\end{equation}
As before, $u>0$ in $\Omega$ and $u_\nu<0$ on $\partial\Omega$, and we write
$u_{\nu\nu}=\nabla^2u(\nu,\nu)$. In this generality the condition
$\operatorname{Ric}^\alpha_w\geq(n+\alpha-1)k\,g$ makes sense for every $k\in\mathbb{R}$ (for
instance, if $w$ is constant and $\operatorname{Ric}\geq(n-1)\kappa\,g$, it holds with
$k=\frac{(n-1)\kappa}{n+\alpha-1}$).

Associated with the solution $u$, we define the \textit{weighted torsion} of $\Omega$ as
\begin{align}\label{torsion}
    T_w(\Omega)=\int_{\Omega} u\,d\mu_w=\int_\Omega|\nabla u|^2\,d\mu_w,
\end{align}
where the second equality follows from
$0=\int_{\Omega}\operatorname{div}_w(u\nabla u)\,d\mu_w=\int_{\Omega}(-u+|\nabla u|^2)\,d\mu_w$.
To derive our geometric estimates, we use the same $P$-function as defined in (\ref{PFunction}).
By applying the maximum principle and the weighted Bochner formula to this $P$-function, we establish a sharp lower bound for the maximum of $|\nabla u|^2$ on the boundary. This is the content of our first theorem.

The following observation shows that, for smooth positive weights, the weighted Bochner
inequality is never an equality near a maximum point of $u$.

\begin{lemma}\label{lem:strict}
Let $u$ be the solution of \eqref{eq:torsion_general} and let $p\in\Omega$ be a maximum point
of $u$. Then there is a neighborhood $U\subset\Omega$ of $p$ where the inequalities in
Lemmas~\ref{Weightedbochner} and~\ref{BochnerP} are strict.
\end{lemma}

\begin{proof}
Since $\nabla u(p)=0$, we have $\frac{\langle\nabla\log w,\nabla u\rangle}{\alpha}(p)=0$, while
$\frac{\Delta_wu}{n+\alpha}(p)=-\frac{1}{n+\alpha}$. Hence the equality condition in
Lemma~\ref{BochnerP}, which is also part of the equality condition in
Lemma~\ref{Weightedbochner}, fails at $p$, and by continuity it fails near $p$.
\end{proof}

\begin{theorem}\label{thm:max}
Let $u$ be the solution of \eqref{eq:torsion_general} and assume that
$\operatorname{Ric}^\alpha_w\geq(n+\alpha-1)k\,g$ in $\Omega$, for some $k\in\mathbb{R}$.
Then
\begin{equation}\label{eq001}
\max_{\partial \Omega}|\nabla u|^2> \dfrac{n+\alpha+2}{n+\alpha}\dfrac{T_w(\Omega)}{Vol_w(\Omega)}+\dfrac{2(n+\alpha-1)k}{Vol_w(\Omega)}\int_{\Omega}u|\nabla u|^2d\mu_w.
\end{equation}
\end{theorem}

\begin{proof}
We start with the integral identity obtained from Green's identity, since $u=0$ and
$P=u_\nu^2$ on $\partial\Omega$:
\begin{equation*}
\int_{\Omega}(u\Delta_wP-P\Delta_wu)d\mu_w=-\int_{\partial\Omega}u_{\nu}^3d\sigma_w.
\end{equation*}
From the weighted Bochner formula in Lemma \ref{BochnerP} and the curvature assumption, we have the differential inequality
\[
\Delta_wP\geq 2\left[|\mathring{\nabla}^2u|^2+Ric_w^{\alpha}(\nabla u,\nabla u)\right]\geq 2Ric_w^{\alpha}(\nabla u,\nabla u)\geq 2(n+\alpha-1)k|\nabla u|^2,
\]
where the first inequality is strict near a maximum point of $u$, by
Lemma~\ref{lem:strict}.
Let $m=\max_{\partial \Omega}|\nabla u|$. Using that $u_{\nu}=-|\nabla u|<0$ (by the Hopf lemma \cite{GT}),
we have $-\int_{\partial\Omega}u_\nu^3\,d\sigma_w\leq-m^2\int_{\partial\Omega}u_\nu\,d\sigma_w$. Multiplying the
inequality above by $u>0$ and integrating, we deduce
\[
2(n+\alpha-1)k\int_{\Omega}u|\nabla u|^2d\mu_w+\int_{\Omega}Pd\mu_w<-m^2\int_{\partial\Omega}u_{\nu}d\sigma_w=m^2Vol_w(\Omega).
\]
On the other hand, by \eqref{torsion}, we obtain
\begin{eqnarray*}
\int_{\Omega}P\,d\mu_w&=&\int_{\Omega}\left(|\nabla u|^2+\dfrac{2}{n+\alpha}u\right)d\mu_w\\
&=&\int_{\Omega}\left(u+\dfrac{2}{n+\alpha}u\right)d\mu_w=\left[\dfrac{n+\alpha+2}{n+\alpha}\right]T_w(\Omega).
\end{eqnarray*}
Therefore, substituting this back, we obtain
\begin{equation}\label{eq002}
m^2> \left[\dfrac{n+\alpha+2}{n+\alpha}\right]\dfrac{T_w(\Omega)}{Vol_w(\Omega)}+\dfrac{2(n+\alpha-1)k}{Vol_w(\Omega)}\int_{\Omega}u|\nabla u|^2d\mu_w.
\end{equation}
\end{proof}

The next result extends \cite[Theorem~1.6]{wang2020} to the drift Laplacian. It provides
upper bounds for the minimum of the second normal derivative of $u$ on $\partial\Omega$ and for
its boundary integral.

\begin{theorem}\label{thm:wangxia} Let $u$ be the solution of \eqref{eq:torsion_general}
and assume that $Ric^\alpha_w\geq (n+\alpha-1)kg$ in $\Omega$. Then:
\begin{enumerate}[i)]
\item
\[
\min_{x\in\partial \Omega} u_{\nu\nu}< -\dfrac{(n+\alpha-1)kT_w(\Omega)}{Vol_{w}(\Omega)}-\dfrac{1}{n+\alpha}.
\]
\item If $H_w\geq 0$ on $\partial\Omega$ and $H_w\not\equiv0$, then
\[
\int_{\partial\Omega}u_{\nu\nu}d\sigma_w<(n+\alpha-1)^{1/2}\left[\dfrac{Vol_w(\Omega)}{n+\alpha}-kT_w(\Omega)\right]^{1/2}\left[\int_{\partial\Omega}H_w\,d\sigma_w\right]^{1/2}-Area_w(\partial\Omega).
\]
\end{enumerate}
\end{theorem}


\begin{proof}
To prove item (i), we first note that $u_{\nu}<0$. By the weighted Bochner formula in Lemma \ref{Weightedbochner}, we have
\begin{eqnarray*}
\dfrac{1}{2}\Delta_w|\nabla u|^2&\geq&\dfrac{(\Delta_wu)^2}{n+\alpha}+\langle\nabla u,\nabla\Delta_wu\rangle+Ric^{\alpha}_w(\nabla u, \nabla u)\\
&\geq& \dfrac{1}{n+\alpha}+(n+\alpha-1)k|\nabla u|^2,
\end{eqnarray*}
where the first inequality is strict near a maximum point of $u$ (Lemma~\ref{lem:strict}).
Integrating the above inequality over $\Omega$ with respect to the weighted measure, and
using \eqref{torsion}, we get
\begin{equation}\label{eq003}
(n+\alpha-1)kT_w(\Omega)+\dfrac{Vol_w(\Omega)}{n+\alpha}<\dfrac{1}{2}\int_{\partial\Omega}\langle\nabla|\nabla u|^2,\nu\rangle d\sigma_w.
\end{equation}
Now, working with the right-hand side of (\ref{eq003}), we obtain
\begin{eqnarray*}
\dfrac{1}{2}\int_{\partial \Omega}\langle\nabla|\nabla u|^2,\nu\rangle d\sigma_w&=&\int_{\partial\Omega}\nabla^2u(\nabla u, \nu)d\sigma_w\\
&=&\int_{\partial\Omega}u_{\nu}u_{\nu\nu}d\sigma_w\\
&\leq&l\int_{\partial\Omega}u_{\nu}d\sigma_w\\
&=&l\int_\Omega\Delta_w ud\mu_w=-lVol_w(\Omega),
\end{eqnarray*}
where we used that $\nabla u=u_{\nu}\nu$, since $u=0$ on $\partial \Omega$, and $l=\displaystyle\min_{x\in\partial \Omega} u_{\nu\nu}$.
Substituting this inequality into \eqref{eq003}, we arrive at the desired inequality in item (i).

To prove item (ii), we use the decomposition of the weighted Laplacian on the boundary
(Lemma~\ref{lem:laplacians}):
\[
\Delta_wu=\bar{\Delta}_wu+H_wu_{\nu}+u_{\nu\nu}.
\]
Integrating over $\partial\Omega$ and using that $\bar{\Delta}_wu=0$ (since $u=0$ on $\partial\Omega$), we have
\begin{equation}\label{eq004}
Area_w(\partial \Omega)+\int_{\partial\Omega}u_{\nu\nu}d\sigma_w=-\int_{\partial\Omega}H_wu_{\nu}d\sigma_w.
\end{equation}
Using the weighted Reilly inequality (Proposition \ref{Reilly}), we deduce
\begin{eqnarray}\label{eq005}
\int_{\partial\Omega}H_wu_{\nu}^2d\sigma_w&<&\dfrac{n+\alpha-1}{n+\alpha}\int_{\Omega}(\Delta_wu)^2d\mu_w-\int_{\Omega}Ric^{\alpha}_w(\nabla u,\nabla u)d\mu_w\nonumber\\
&\leq&\dfrac{n+\alpha-1}{n+\alpha}Vol_w(\Omega)-(n+\alpha-1)k\int_{\Omega}|\nabla u|^2d\mu_w\nonumber\\
&=&\dfrac{n+\alpha-1}{n+\alpha}Vol_w(\Omega)-(n+\alpha-1)kT_w(\Omega),
\end{eqnarray}
where the first inequality is strict because, by Lemma~\ref{lem:strict}, equality cannot
hold in Proposition~\ref{Reilly}.
Since $H_w\geq 0$ and $u_{\nu}<0$, we can use the Cauchy--Schwarz   
inequality:
\[
-\int_{\partial \Omega}H_wu_{\nu}d\sigma_w\leq\left[\int_{\partial \Omega}H_wu_{\nu}^2d\sigma_w\right]^{1/2}\left[\int_{\partial \Omega}H_{w}d\sigma_w\right]^{1/2}.
\]
Substituting into \eqref{eq004} and using \eqref{eq005}, together with
$\int_{\partial\Omega}H_w\,d\sigma_w>0$, we deduce
\begin{eqnarray*}
\int_{\partial \Omega}u_{\nu\nu}d\sigma_w&\leq&\left(\int_{\partial\Omega}H_wu_{\nu}^2d\sigma_w\right)^{1/2}\left(\int_{\partial\Omega}H_wd\sigma_w\right)^{1/2}-Area_w(\partial \Omega)\\
&<&(n+\alpha-1)^{1/2}\left[\dfrac{Vol_w(\Omega)}{n+\alpha}-kT_w(\Omega)\right]^{1/2}\left[\int_{\partial\Omega}H_wd\sigma_w\right]^{1/2}-Area_w(\partial \Omega).
\end{eqnarray*}
\end{proof}

\begin{remark}\label{rem:strict}
In the unweighted case \cite{wang2020}, equality in these estimates holds if and only if $k=0$
and $\Omega$ is isometric to a Euclidean ball. The strictness in Theorems~\ref{thm:max}
and~\ref{thm:wangxia} shows that this rigidity is specific to the unweighted case: for a
smooth positive weight, the inequality in Lemma~\ref{BochnerP} cannot be an equality at the
maximum point of $u$, where $\nabla u$ vanishes but $\Delta_wu$ does not. Equality requires the
weight to degenerate there, and this is exactly what happens for homogeneous weights in
$\mathbb{R}^n$, to which we now return. 
\end{remark}

\section{A weighted general identity and the reverse Serrin problem}

In this section we study the weighted reverse Serrin problem described in the
Introduction \cite{Magnanini2024}: we start from a solution of the Neumann problem
\begin{equation}\label{eq:reverse_neumann}
    \Delta_w u=-1\quad\text{in }\Omega,
    \qquad
    u_\nu=-r\quad\text{on }\partial\Omega ,
    \qquad\text{where}\qquad
   r=\frac{Vol_w(\Omega)}{Area_w(\partial\Omega)}
\end{equation}
is forced by the weighted divergence theorem; such a solution exists for every domain
and is unique up to an additive constant. We look for conditions on the boundary values of
$u$, weaker than the Dirichlet one, which force $\Omega$ to be a ball.
Following \cite{Magnanini2024},
we first establish a general integral identity for two solutions of the
weighted Poisson equation, valid on any weighted manifold, and then specialize it to
\eqref{eq:reverse_neumann} in $\mathbb{R}^n$.

\begin{proposition}[General weighted identity]\label{thm:general_weighted}
Let $\Omega$ be a bounded domain with smooth boundary in a weighted Riemannian manifold
$(M^n,\langle\cdot,\cdot\rangle,w\,d\mathrm{vol})$, and let $U,V\in C^2(\overline{\Omega})\cap C^3(\Omega)$ satisfy
\begin{equation*}
    \Delta_w U=\Delta_w V=\Lambda\quad\text{in }\Omega,
\end{equation*}
for some constant $\Lambda\in\mathbb{R}$. Set
\begin{equation*}
    \overline{U}=\max_{\overline{\Omega}}U
    \qquad\text{and}\qquad
    \tilde P=\frac12|\nabla U|^2+(\overline{U}-U)\quad\text{in }\overline{\Omega}.
\end{equation*}
Then
\begin{align}\label{eq:general_identity}
    \int_{\Omega}(\overline{U}-U)\,\Delta_w\tilde P\,d\mu_w
    &+\int_{\Omega}\langle(I-\nabla^2V)\nabla U,\nabla U\rangle\,d\mu_w \nonumber\\
    &=\int_{\partial\Omega}(\overline{U}-U)\,\nabla^2U(\nabla U,\nu)\,d\sigma_w
      +\frac12\int_{\partial\Omega}|\nabla U|^2(U_\nu+V_\nu)\,d\sigma_w \nonumber\\
    &\quad-\int_{\partial\Omega}\langle\nabla V,\nabla U\rangle U_\nu\,d\sigma_w
      -\int_{\partial\Omega}(\overline{U}-U)\,U_\nu\,d\sigma_w \nonumber\\
    &\quad+\Lambda\int_{\partial\Omega}(\overline{U}-U)(U_\nu-V_\nu)\,d\sigma_w .
\end{align}
Here \(I\) denotes the identity endomorphism of \(TM\), and \(\nabla^2V\) is identified with an endomorphism via the Riemannian metric.
\end{proposition}

\begin{remark}
The proof below uses only that $\overline{U}$ is a constant, so
\eqref{eq:general_identity} holds for any $\overline{U}\in\mathbb{R}$.
\end{remark}

\begin{proof}
Since $\overline U$ is constant, $\nabla(\overline U-U)=-\nabla U$. Consider the
vector fields
\begin{align*}
    X&=\tilde P\,\nabla U+(\overline U-U)\,\nabla\tilde P+\frac12|\nabla U|^2\nabla V
       -\langle\nabla V,\nabla U\rangle\nabla U,\\
    Y&=(\Lambda-1)(\overline U-U)\,\nabla U-\Lambda(\overline U-U)\,\nabla V .
\end{align*}
Using the product rule for $\operatorname{div}_w$, we compute each term:
\begin{align*}
    \operatorname{div}_w(\tilde P\,\nabla U)
      &=\langle\nabla\tilde P,\nabla U\rangle+\tilde P\,\Delta_wU,\\
    \operatorname{div}_w\big((\overline U-U)\,\nabla\tilde P\big)
      &=-\langle\nabla U,\nabla\tilde P\rangle+(\overline U-U)\,\Delta_w\tilde P,\\
    \operatorname{div}_w\Big(\frac12|\nabla U|^2\nabla V\Big)
      &=\nabla^2U(\nabla U,\nabla V)+\frac12|\nabla U|^2\Delta_wV,\\
    \operatorname{div}_w\big(-\langle\nabla V,\nabla U\rangle\nabla U\big)
      &=-\nabla^2V(\nabla U,\nabla U)-\nabla^2U(\nabla V,\nabla U)
        -\langle\nabla V,\nabla U\rangle\Delta_wU,\\
    \operatorname{div}_w\big((\Lambda-1)(\overline U-U)\,\nabla U\big)
      &=-(\Lambda-1)|\nabla U|^2+(\Lambda-1)(\overline U-U)\,\Delta_wU,\\
    \operatorname{div}_w\big(-\Lambda(\overline U-U)\,\nabla V\big)
      &=\Lambda\langle\nabla U,\nabla V\rangle-\Lambda(\overline U-U)\,\Delta_wV .
\end{align*}
Adding these identities and using the fact that $\Delta_wU=\Delta_wV=\Lambda$, the terms
$\pm\langle\nabla U,\nabla\tilde P\rangle$, $\pm\nabla^2U(\nabla U,\nabla V)$ and
$\pm \Lambda\langle\nabla U,\nabla V\rangle$ cancel, yielding
\begin{align*}
    \operatorname{div}_w(X+Y)
    &=(\overline U-U)\,\Delta_w\tilde P-\nabla^2V(\nabla U,\nabla U)
      +\Lambda\tilde P+\frac \Lambda2|\nabla U|^2-(\Lambda-1)|\nabla U|^2\\
    &\quad+(\Lambda-1)\Lambda(\overline U-U)-\Lambda^2(\overline U-U)\\
    &=(\overline U-U)\,\Delta_w\tilde P-\nabla^2V(\nabla U,\nabla U)
      +\Lambda\Big(\tilde P-\frac12|\nabla U|^2-(\overline U-U)\Big)+|\nabla U|^2\\
    &=(\overline{U}-U)\,\Delta_w\tilde P+\langle(I-\nabla^2V)\nabla U,\nabla U\rangle,
\end{align*}
where in the last step we used $\tilde P=\frac12|\nabla U|^2+(\overline U-U)$.

We now integrate over $\Omega$ with respect to $d\mu_w$ and apply the weighted
divergence theorem. On $\partial\Omega$ we have
$\tilde P_\nu=\nabla^2U(\nabla U,\nu)-U_\nu$, hence
\begin{align*}
    \langle X+Y,\nu\rangle
    &=\Big(\frac12|\nabla U|^2+(\overline U-U)\Big)U_\nu
      +(\overline U-U)\big(\nabla^2U(\nabla U,\nu)-U_\nu\big)\\
    &\quad+\frac12|\nabla U|^2V_\nu-\langle\nabla V,\nabla U\rangle U_\nu
      +(\Lambda-1)(\overline U-U)\,U_\nu-\Lambda(\overline U-U)\,V_\nu\\
    &=(\overline U-U)\,\nabla^2U(\nabla U,\nu)+\frac12|\nabla U|^2(U_\nu+V_\nu)
      -\langle\nabla V,\nabla U\rangle U_\nu\\
    &\quad-(\overline U-U)\,U_\nu+\Lambda(\overline U-U)(U_\nu-V_\nu),
\end{align*}
which gives \eqref{eq:general_identity}.
\end{proof}

From now on $\Omega\subset\mathbb{R}^n$ and $w$ satisfy the standing assumptions of Section~2;
in particular $w$ is homogeneous of degree $\alpha>0$, so that
$\langle\nabla\log w,x\rangle=\alpha$. In this case the
quadratic polynomial
\begin{equation*}
    V(x)=\frac{|x|^2}{2}
    \qquad\text{satisfies}\qquad
    \nabla^2V=I
    \quad\text{and}\quad
    \Delta_wV=n+\langle\nabla\log w,x\rangle=n+\alpha ,
\end{equation*}
and it will play the role of the model solution. Accordingly, we normalize
$U$ so that $\Delta_wU=n+\alpha$; this corresponds to the problem
$\Delta_wu=-1$ of Section~3 through the change of variables
$U=-(n+\alpha)u$.

Taking $\Lambda=n+\alpha$ and this choice of $V$ in
Proposition~\ref{thm:general_weighted}, the second integral on the left-hand
side of \eqref{eq:general_identity} vanishes, because $\nabla^2V=I$. Moreover,
the last two boundary integrals in \eqref{eq:general_identity} combine into
\begin{equation*}
    -(\overline U-U)\,U_\nu+(n+\alpha)(\overline U-U)(U_\nu-V_\nu)
    =(\overline U-U)\big[(n+\alpha-1)U_\nu-(n+\alpha)V_\nu\big].
\end{equation*}
Finally, splitting $\nabla U=\bar\nabla U+U_\nu\nu$ and
$\nabla V=\bar\nabla V+V_\nu\nu$ on $\partial\Omega$, we have, pointwise on
$\partial\Omega$,
\begin{equation*}
    \frac12|\nabla U|^2(U_\nu+V_\nu)-\langle\nabla V,\nabla U\rangle U_\nu
    =\frac12U_\nu^2(U_\nu-V_\nu)+\frac12|\bar\nabla U|^2(U_\nu+V_\nu)
      -\langle\bar\nabla V,\bar\nabla U\rangle U_\nu .
\end{equation*}
This yields the following form of the identity, in which the normal and
tangential parts of $\nabla U$ appear separately.

\begin{corollary}\label{cor:alternative}
Let $V=|x|^2/2$, and let
$U\in C^2(\overline\Omega)$ satisfy $\Delta_wU=n+\alpha$ in $\Omega$. Then
\begin{align}\label{teo 5*}
    \int_{\Omega}(\overline{U}-U)\,\Delta_w\tilde P\,d\mu_w
    &=\frac12\int_{\partial\Omega}U_\nu^2(U_\nu-V_\nu)\,d\sigma_w
      +\frac12\int_{\partial\Omega}|\bar\nabla U|^2(U_\nu+V_\nu)\,d\sigma_w \nonumber\\
    &\quad-\int_{\partial\Omega}\langle\bar\nabla V,\bar\nabla U\rangle U_\nu\,d\sigma_w
      +\int_{\partial\Omega}(\overline{U}-U)\,\nabla^2U(\nabla U,\nu)\,d\sigma_w \nonumber\\
    &\quad+\int_{\partial\Omega}(\overline{U}-U)\big[(n+\alpha-1)U_\nu-(n+\alpha)V_\nu\big]\,d\sigma_w .
\end{align}
\end{corollary}

We now turn to the reverse Serrin problem, that is, we assume that the
Neumann datum is constant:
\begin{equation}\label{eq:neumann_U}
    U_\nu=R\quad\text{on }\partial\Omega,
    \qquad\text{where}\qquad
    R=(n+\alpha)\,\frac{Vol_w(\Omega)}{Area_w(\partial\Omega)}
\end{equation}
is forced by the weighted divergence theorem. In terms of $u$, this is exactly
the condition $u_\nu=-r$ of \eqref{eq:reverse_neumann}, with $R=(n+\alpha)\,r$. One key
step is to rewrite the term $\nabla^2U(\nabla U,\nu)$, which involves second
derivatives of $U$ in the normal direction. This is done in the next lemma,
which is the weighted counterpart of \cite[Lemma~2.1]{Magnanini2024}.

\begin{lemma}\label{lema:weighted_2.1}
Let $U\in C^2(\overline\Omega)$ satisfy $\Delta_wU=n+\alpha$ in $\Omega$ and
$U_\nu=R$ on $\partial\Omega$, for some constant $R$. Then, on $\partial\Omega$,
\begin{equation}\label{eq:lemma21}
    \nabla^2U(\nabla U,\nu)
    =-h(\bar\nabla U,\bar\nabla U)
     +U_\nu\big[(n+\alpha)-\bar\Delta_wU-H_wU_\nu\big],
\end{equation}
where $h$ and $H_w$ are the second fundamental form and the weighted mean
curvature of $\partial\Omega$ with respect to $\nu$.
\end{lemma}

\begin{proof}
Extend $\nu$ to a unit vector field of class $C^1$ in a neighborhood of
$\partial\Omega$. Since $U_\nu=\langle\nabla U,\nu\rangle=R$ is constant along
$\partial\Omega$, its tangential gradient vanishes there, so its full gradient
is normal to the boundary:
\begin{equation}\label{derivada gradu}
    \nabla\langle\nabla U,\nu\rangle
    =\big\langle\nabla\langle\nabla U,\nu\rangle,\nu\big\rangle\,\nu
    \quad\text{on }\partial\Omega .
\end{equation}
Moreover, differentiating $\langle\nu,\nu\rangle=1$ gives
$\langle\nabla_X\nu,\nu\rangle=0$ for every vector field $X$. In particular,
$\nabla_X\nu$ is tangent to $\partial\Omega$, and for $X$ tangent we have
$\langle\nabla_X\nu,\bar\nabla U\rangle=h(X,\bar\nabla U)$.

We compute the derivative of $\langle\nabla U,\nu\rangle$ in the direction of
$\nabla U$. On $\partial\Omega$,
\begin{align*}
    \big\langle\nabla\langle\nabla U,\nu\rangle,\nabla U\big\rangle
    &=\nabla^2U(\nabla U,\nu)+\langle\nabla U,\nabla_{\nabla U}\nu\rangle\\
    &=\nabla^2U(\nabla U,\nu)+\langle\bar\nabla U,\nabla_{\nabla U}\nu\rangle
      +U_\nu\langle\nu,\nabla_{\nabla U}\nu\rangle\\
    &=\nabla^2U(\nabla U,\nu)+\langle\bar\nabla U,\nabla_{\nabla U}\nu\rangle\\
    &=\nabla^2U(\nabla U,\nu)+\langle\bar\nabla U,\nabla_{\bar\nabla U}\nu\rangle
      +U_\nu\langle\bar\nabla U,\nabla_\nu\nu\rangle\\
    &=\nabla^2U(\nabla U,\nu)+h(\bar\nabla U,\bar\nabla U)
      +U_\nu\langle\bar\nabla U,\nabla_\nu\nu\rangle .
\end{align*}
Here we decomposed $\nabla U=\bar\nabla U+U_\nu\nu$ and used
$\langle\nu,\nabla_X\nu\rangle=0$ and $\langle\nabla_X\nu,Y\rangle=h(X,Y)$ for $X,Y$ tangent to
$\partial\Omega$.
On the other hand, by \eqref{derivada gradu}, Lemma~\ref{lem:laplacians} and
$\Delta_wU=n+\alpha$,
\begin{align*}
    \big\langle\nabla\langle\nabla U,\nu\rangle,\nabla U\big\rangle
    &=\big\langle\nabla\langle\nabla U,\nu\rangle,\nu\big\rangle\,U_\nu
     =\nu\big(\langle\nabla U,\nu\rangle\big)\,U_\nu\\
    &=\big(\nabla^2U(\nu,\nu)+\langle\nabla U,\nabla_\nu\nu\rangle\big)U_\nu\\
    &=\big((n+\alpha)-\bar\Delta_wU-H_wU_\nu+\langle\bar\nabla U,\nabla_\nu\nu\rangle\big)U_\nu,
\end{align*}
where we used
$\langle\nabla U,\nabla_\nu\nu\rangle
=U_\nu\langle\nu,\nabla_\nu\nu\rangle+\langle\bar\nabla U,\nabla_\nu\nu\rangle
=\langle\bar\nabla U,\nabla_\nu\nu\rangle$.
Comparing the two expressions, the terms
$U_\nu\langle\bar\nabla U,\nabla_\nu\nu\rangle$ cancel, and we obtain
\eqref{eq:lemma21}. In particular, the result does not depend on the chosen
extension of $\nu$.
\end{proof}

Next we establish the weighted version of \cite[Corollary~3.3]{Magnanini2024},
which simplifies \eqref{teo 5*} when the Neumann datum is constant.

\begin{proposition}\label{prop:weighted_3.3}
Let $V=|x|^2/2$, and let $U\in C^2(\overline\Omega)$ satisfy
\begin{equation*}
    \Delta_wU=n+\alpha\quad\text{in }\Omega,
    \qquad
    U_\nu=R\quad\text{on }\partial\Omega .
\end{equation*}
Then
\begin{align}\label{eq:weighted_3.3}
    \int_\Omega(\overline U-U)\,\Delta_w\tilde P\,d\mu_w
    &=\frac12\int_{\partial\Omega}|\bar\nabla U|^2(V_\nu-R)\,d\sigma_w
      -\int_{\partial\Omega}(\overline U-U)\,h(\bar\nabla U,\bar\nabla U)\,d\sigma_w \nonumber\\
    &\quad+\int_{\partial\Omega}(\overline U-U)\big[RH_w-(n+\alpha)\big](V_\nu-R)\,d\sigma_w .
\end{align}
\end{proposition}

\begin{proof}
We examine the right-hand side of \eqref{teo 5*}, which we denote by $\mathrm{I}+\mathrm{II}+\mathrm{III}+\mathrm{IV}+\mathrm{V}$, following the order in which the terms appear.

\smallskip\noindent
First, we evaluate term $\mathrm{I}$. Since $U_\nu=R$ and $\Delta_wU=\Delta_wV$, the divergence theorem yields
\begin{equation*}
    \mathrm{I}=\frac{R^2}{2}\int_{\partial\Omega}(U_\nu-V_\nu)\,d\sigma_w
    =\frac{R^2}{2}\int_\Omega(\Delta_wU-\Delta_wV)\,d\mu_w=0 .
\end{equation*}

\smallskip\noindent
Next, we compute $\mathrm{II}+\mathrm{III}$. Since $\overline U$ is constant, we have $\bar\nabla(\overline U-U)=-\bar\nabla U$. Integrating by parts on $\partial\Omega$, we obtain
\begin{align*}
    \mathrm{III}&=-R\int_{\partial\Omega}\langle\bar\nabla V,\bar\nabla U\rangle\,d\sigma_w
    =R\int_{\partial\Omega}\langle\bar\nabla V,\bar\nabla(\overline U-U)\rangle\,d\sigma_w\\
    &=-R\int_{\partial\Omega}(\overline U-U)\,\bar\Delta_wV\,d\sigma_w .
\end{align*}
By Lemma~\ref{lem:laplacians}, together with $\Delta_wV=n+\alpha$ and $\nabla^2V(\nu,\nu)=1$, it follows that
\begin{equation*}
    \bar\Delta_wV=\Delta_wV-\nabla^2V(\nu,\nu)-H_wV_\nu=(n+\alpha-1)-H_wV_\nu .
\end{equation*}
Therefore,
\begin{equation*}
    \mathrm{II}+\mathrm{III}
    =\frac12\int_{\partial\Omega}|\bar\nabla U|^2(R+V_\nu)\,d\sigma_w
     -R\int_{\partial\Omega}(\overline U-U)\big[(n+\alpha-1)-H_wV_\nu\big]\,d\sigma_w.
\end{equation*}

\smallskip\noindent
For term $\mathrm{IV}$, combining integration by parts on $\partial\Omega$, namely
$\int_{\partial\Omega}(\overline U-U)\,\bar\Delta_wU\,d\sigma_w=\int_{\partial\Omega}|\bar\nabla U|^2\,d\sigma_w$,
with Lemma~\ref{lema:weighted_2.1} gives
\begin{align*}
    \mathrm{IV}
    &=-\int_{\partial\Omega}(\overline U-U)\,h(\bar\nabla U,\bar\nabla U)\,d\sigma_w
     -R\int_{\partial\Omega}|\bar\nabla U|^2\,d\sigma_w\\
    &\quad-R^2\int_{\partial\Omega}(\overline U-U)H_w\,d\sigma_w
     +(n+\alpha)R\int_{\partial\Omega}(\overline U-U)\,d\sigma_w.
\end{align*}

\smallskip\noindent
Finally, since $U_\nu=R$, term $\mathrm{V}$ becomes
\begin{equation*}
    \mathrm{V}=\int_{\partial\Omega}(\overline U-U)\big[(n+\alpha-1)R-(n+\alpha)V_\nu\big]\,d\sigma_w.
\end{equation*}

\smallskip\noindent
We now group the resulting integrals. Collecting the terms involving $|\bar\nabla U|^2$ gives
\begin{equation*}
    \frac12\int_{\partial\Omega}|\bar\nabla U|^2(R+V_\nu)\,d\sigma_w
    -R\int_{\partial\Omega}|\bar\nabla U|^2\,d\sigma_w
    =\frac12\int_{\partial\Omega}|\bar\nabla U|^2(V_\nu-R)\,d\sigma_w,
\end{equation*}
while the terms involving $H_w$ combine to yield
\begin{equation*}
    R\int_{\partial\Omega}(\overline U-U)H_wV_\nu\,d\sigma_w
    -R^2\int_{\partial\Omega}(\overline U-U)H_w\,d\sigma_w
    =R\int_{\partial\Omega}(\overline U-U)H_w(V_\nu-R)\,d\sigma_w.
\end{equation*}
In the remaining terms, the coefficients of $R\int_{\partial\Omega}(\overline U-U)\,d\sigma_w$ arising from $\mathrm{II}+\mathrm{III}$, $\mathrm{IV}$, and $\mathrm{V}$ sum to $-(n+\alpha-1)+(n+\alpha)+(n+\alpha-1)=n+\alpha$. Hence, these remaining integrals simplify to
\begin{align*}
    &(n+\alpha)R\int_{\partial\Omega}(\overline U-U)\,d\sigma_w
    -(n+\alpha)\int_{\partial\Omega}(\overline U-U)V_\nu\,d\sigma_w\\
    &\qquad=-(n+\alpha)\int_{\partial\Omega}(\overline U-U)(V_\nu-R)\,d\sigma_w.
\end{align*}
Summing these three expressions together with the term containing $h$ yields \eqref{eq:weighted_3.3}.
\end{proof}

We now translate Proposition~\ref{prop:weighted_3.3} in terms of the solution
$u$ of the reverse Serrin problem \blue{\eqref{eq:reverse_neumann}}.
Setting $U=-(n+\alpha)u$, we have $\Delta_wU=n+\alpha$ in $\Omega$ and
$U_\nu=R$ on $\partial\Omega$, with
\begin{equation*}
    R=(n+\alpha)\,r .
\end{equation*}
Let $u_*=\min_{\overline\Omega}u$, so that $\overline U-U=(n+\alpha)(u-u_*)$,
and let $P$ be the $P$-function \eqref{PFunction}.
Then $\tilde P=\frac{(n+\alpha)^2}{2}\,P-(n+\alpha)u_*$, and hence
$\Delta_w\tilde P=\frac{(n+\alpha)^2}{2}\,\Delta_wP$. Moreover,
$V=-(n+\alpha)q$, where $q(x)=-\frac{|x|^2}{2(n+\alpha)}$ up to an additive
constant, which plays no role. Since $q_\nu=-\langle x,\nu\rangle/(n+\alpha)$,
we get
\begin{equation*}
    V_\nu-R=\langle x,\nu\rangle-(n+\alpha)r=(n+\alpha)(u_\nu-q_\nu),
    \quad
    RH_w-(n+\alpha)=-(n+\alpha)\,(1-r H_w)
    \text{ on }\partial\Omega .
\end{equation*}

\begin{corollary}\label{cor:weighted_3.3_u}
Let $u\in C^2(\overline\Omega)$
be a solution of \eqref{eq:reverse_neumann}. Then
\begin{align}\label{eq:weighted_3.3_u}
    \int_\Omega(u-u_*)\,\Delta_wP\,d\mu_w
    &=\int_{\partial\Omega}|\bar\nabla u|^2(u_\nu-q_\nu)\,d\sigma_w
      -2\int_{\partial\Omega}(u-u_*)\,h(\bar\nabla u,\bar\nabla u)\,d\sigma_w \nonumber\\
    &\quad-2\int_{\partial\Omega}(u-u_*)\,(1-r H_w)(u_\nu-q_\nu)\,d\sigma_w.
\end{align}
\end{corollary}

\begin{proof}
We begin by substituting the relations established above into \eqref{eq:weighted_3.3}. Since
$|\bar\nabla U|^2=(n+\alpha)^2|\bar\nabla u|^2$,
$h(\bar\nabla U,\bar\nabla U)=(n+\alpha)^2h(\bar\nabla u,\bar\nabla u)$ and
$\Delta_w\tilde P=\frac{(n+\alpha)^2}{2}\Delta_wP$, we have
\begin{gather*}
    (\overline U-U)\,\Delta_w\tilde P=\frac{(n+\alpha)^3}{2}(u-u_*)\,\Delta_wP,\qquad
    \frac12|\bar\nabla U|^2(V_\nu-R)=\frac{(n+\alpha)^3}{2}|\bar\nabla u|^2(u_\nu-q_\nu),\\
    (\overline U-U)\,h(\bar\nabla U,\bar\nabla U)=(n+\alpha)^3(u-u_*)\,h(\bar\nabla u,\bar\nabla u),\\
    (\overline U-U)\big[RH_w-(n+\alpha)\big](V_\nu-R)
    =-(n+\alpha)^3(u-u_*)(1-r H_w)(u_\nu-q_\nu).
\end{gather*}
Dividing by $(n+\alpha)^3/2$ gives \eqref{eq:weighted_3.3_u}.
\end{proof}


As in \cite[Lemma~3.4]{Magnanini2024}, the key point in the rigidity argument is that the integral identity, together with the curvature assumption, forces the traceless Hessian of $u$ to vanish. In the weighted setting, the vanishing of the Bergström term in Lemma~\ref{BochnerP} also plays a role: it fixes the dilation constant of the quadratic polynomial.

\begin{lemma}\label{lem:rigidity}
Let $u\in C^3(\overline\Omega)$
be a solution of \eqref{eq:reverse_neumann}.
Assume that
\[
\operatorname{Ric}^\alpha_w(\nabla u,\nabla u)
\geq 0
\quad\text{in }\Omega\setminus\{0\}.
\]
Let $P$ be given by \eqref{PFunction}.
Then
\[
\Delta_wP\geq0
\quad\text{in }\Omega\setminus\{0\}.
\]
If, moreover, $\Delta_wP\equiv0$ in
$\Omega\setminus\{0\}$, then
$\Omega=B_R(0)$, where $R=(n+\alpha)r$, and
\[
u(x)=u_*+\frac{R^2-|x|^2}{2(n+\alpha)}.
\]
\end{lemma}

\begin{proof}
By Lemma~\ref{BochnerP} and the curvature
assumption, we have
\[
\Delta_wP
\geq 2|\mathring{\nabla}^2u|^2
\geq0
\quad\text{in }\Omega\setminus\{0\}.
\]
If $\Delta_wP\equiv0$.
Then $\mathring{\nabla}^2u=0$, and equality
holds in Lemma~\ref{BochnerP}. Hence
\[
\frac{\Delta u}{n}
=
\frac{\langle\nabla\log w,\nabla u\rangle}
{\alpha}
=
\frac{\Delta_wu}{n+\alpha}
=
-\frac1{n+\alpha}.
\]
Consequently,
\[
\nabla^2u=-\frac1{n+\alpha}g.
\]
Since $\Omega$ is connected, there exist
$A\in\mathbb R$ and $x_0\in\mathbb R^n$
such that
\[
u(x)=A-\frac{|x-x_0|^2}{2(n+\alpha)}.
\]

The equality condition also yields
\[
\langle\nabla\log w,x-x_0\rangle=\alpha.
\]
By Euler's identity,
$\langle\nabla\log w,x\rangle=\alpha$,
and therefore
\[
\langle\nabla w,x_0\rangle=0
\quad\text{in }\Omega\setminus\{0\}.
\]

The Neumann boundary condition gives
\[
\langle x-x_0,\nu\rangle
=(n+\alpha)r=R>0
\quad\text{on }\partial\Omega.
\]
This implies that $x_0\in\Omega$.
Indeed, if $x_0\notin\Omega$, considering
a boundary point closest to $x_0$
contradicts the preceding identity.

If $x_0\neq0$, evaluating
$\langle\nabla w,x_0\rangle=0$
at $x=x_0$ and using Euler's identity
gives
\[
0=\alpha w(x_0)>0,
\]
a contradiction. Thus $x_0=0$.

The Neumann condition now reads
\[
\langle x,\nu\rangle=R
\quad\text{on }\partial\Omega.
\]
At points of $\partial\Omega$ nearest to
and farthest from the origin, the position
vector is parallel to the outward normal.
Hence the minimum and maximum distances
from the origin to $\partial\Omega$ both
equal $R$. Therefore
\[
\Omega=B_R(0).
\]

Finally, since $u_*$ is attained on
$\partial\Omega$, we obtain
\[
A=u_*+\frac{R^2}{2(n+\alpha)},
\]
which gives the desired expression for $u$.
\end{proof}

\begin{theorem}\label{thm:reverse_serrin}
Let $u\in C^3(\overline\Omega)$ be a solution of \eqref{eq:reverse_neumann} such that
\begin{equation}\label{eq:int_ric_reverse}
    \int_\Omega(u-u_*)\operatorname{Ric}^\alpha_w(\nabla u,\nabla u)\,d\mu_w\geq0 .
\end{equation}
If the
right-hand side of \eqref{eq:weighted_3.3_u} is non-positive, then
$\Omega=B_R(0)$, where $R=(n+\alpha)r$, and
\begin{equation*}
    u(x)=u_*+\frac{R^2-|x|^2}{2(n+\alpha)} .
\end{equation*}
\end{theorem}

\begin{proof}
Since $\Delta_wu=-1<0$, the strong
minimum principle gives $u>u_*$ in
$\Omega$. By Lemma~\ref{BochnerP} and
\eqref{eq:int_ric_reverse},
\[
\begin{aligned}
\int_\Omega (u-u_*)\Delta_wP\,d\mu_w
&\geq
2\int_\Omega (u-u_*)
\Big[
|\mathring{\nabla}^2u|^2
+\operatorname{Ric}^\alpha_w
(\nabla u,\nabla u)
\Big]\,d\mu_w\\
&\geq
2\int_\Omega (u-u_*)
|\mathring{\nabla}^2u|^2\,d\mu_w
\geq0.
\end{aligned}
\]
On the other hand, the left-hand side
is non-positive by
\eqref{eq:weighted_3.3_u} and the
assumption. Hence all these quantities
vanish.

Since $u-u_*>0$ in $\Omega$, we obtain
\[
\mathring{\nabla}^2u=0.
\]
Moreover, the nonnegative function
\[
\Delta_wP-
2\Big[
|\mathring{\nabla}^2u|^2
+\operatorname{Ric}^\alpha_w
(\nabla u,\nabla u)
\Big]
\]
has vanishing integral against the
positive weight $(u-u_*)\,d\mu_w$.
Thus equality holds in
Lemma~\ref{BochnerP} at every point
of $\Omega\setminus\{0\}$.

Consequently,
\[
\frac{\Delta u}{n}
=
\frac{\Delta_wu}{n+\alpha}
=
-\frac{1}{n+\alpha}.
\]
Together with
$\mathring{\nabla}^2u=0$, this yields
\[
\nabla^2u=-\frac{1}{n+\alpha}g.
\]
Since $P=|\nabla u|^2+
\frac{2u}{n+\alpha}$, we compute
\[
\begin{aligned}
\nabla P
&=\nabla(|\nabla u| ^2)
+\frac{2}{n+\alpha}\nabla u\\
&=0.
\end{aligned}
\]
Therefore $P$ is constant in $\Omega$,
and hence
\[
\Delta_wP\equiv0
\quad\text{in }\Omega\setminus\{0\}.
\]

Since equality holds in
Lemma~\ref{BochnerP} and
$\mathring{\nabla}^2u=0$, we obtain
\[
0=\Delta_wP
=2\operatorname{Ric}^\alpha_w
(\nabla u,\nabla u).
\]
In particular,
\[
\operatorname{Ric}^\alpha_w
(\nabla u,\nabla u)=0
\quad\text{in }\Omega\setminus\{0\}.
\]

Thus all the hypotheses of
Lemma~\ref{lem:rigidity} are satisfied.
The conclusion follows.
\end{proof}

The following corollary shows that Theorem~\ref{thm:reverse_serrin} recovers the weighted Serrin rigidity result when the solution is constant on the boundary.

\begin{corollary}\label{cor:reverse_concrete}
Let $u$ be as in Theorem~\ref{thm:reverse_serrin},
satisfying \eqref{eq:int_ric_reverse}.
If $u$ is constant on $\partial\Omega$, then
$\Omega=B_R(0)$, where $R=(n+\alpha)r$, and
\[
u(x)=u_*+\frac{R^2-|x|^2}{2(n+\alpha)}.
\]
\end{corollary}

\begin{proof}
Since $u$ is constant on $\partial\Omega$, we have
$\bar\nabla u=0$ on $\partial\Omega$.
Moreover, as $\Delta_wu=-1<0$, the minimum
principle implies that $u$ attains its minimum
$u_*$ on $\partial\Omega$. Hence $u=u_*$
throughout $\partial\Omega$.

Consequently, every term on the right-hand side
of \eqref{eq:weighted_3.3_u} vanishes, and thus
\[
\int_\Omega (u-u_*)\Delta_wP\,d\mu_w=0.
\]
Therefore, the conclusion follows directly
from Theorem~\ref{thm:reverse_serrin}.
\end{proof}

\section*{Data availability statement}
This manuscript has no associated data.

\section*{Conflict of interest statement}
On behalf of all authors, the corresponding author states that there is no conflict of interest.

\section*{Funding}
M. Andrade was partially supported by the Brazilian National Council for Scientific and Technological Development (CNPq, grants 408834/2023-4, 403869/2024-2, 318922/2025-8 and 400078/2025-2) and FAPITEC/SE/Brazil (grant 019203.01303/2024-1).

A. Freitas was partially supported by the Brazilian National Council for Scientific and Technological Development (CNPq), under Grants No. 308141/2025-3 and No. 406078/2025-4. He also acknowledges support from the Paraíba State Research Foundation (FAPESQ), Brazil, Grant 2453/2026. This study was financed, in part, by the São Paulo Research Foundation (FAPESP), Brazil, Process Number 2026/10743-0.

J. S. Sindeaux was partially supported by the Ceará Foundation for Support of Scientific and Technological Development - FUNCAP/CE, Brazil, under Grant BPI 11/2024: BP6-0241-00436.01.00/25 and by the Brazilian National Council for Scientific and Technological Development (CNPq), Brazil, Grant 406078/2025-4. She also acknowledges support from the Paraíba State Research Foundation (FAPESQ), Brazil, Grant 2453/2026.

\section*{Acknowledgments}
A. Freitas wishes to thank IME/USP and IMECC/UNICAMP for the stimulating scientific environment and support provided during the development of part of this work. In particular, he is grateful to Paolo Piccione and Lino Grama for their warm hospitality and constant encouragement.

\end{document}